\documentstyle[12pt]{article}
\newcommand{\mysection}[1]{\section{#1}\setcounter{equation}{0}}

\title{\bf Iterated logarithm approximations to the distribution of the
largest prime divisor}

\author{{\bf Arie Leizarowitz}
\\{Department of Mathematics}\\{Technion, Haifa 32000}\\
{Israel}}

\date{March 2009}
\begin{document}%

\maketitle

\begin{center} {\bf Abstract} \end{center}
\small
The paper is concerned with estimating the number of integers smaller than $x$ whose largest
prime divisor is smaller than $y$, denoted $\psi (x,y)$. Much of the related literature is
concerned with approximating $\psi (x,y)$ by Dickman's function $\rho (u)$, where
$u=\ln x/\ln y$. A typical such result is that
$$
\psi (x,y)=x\rho (u)(1+o(1)) \eqno (1)
$$ in a certain
domain of the parameters $x$ and $y$.

In this paper a different type of approximation of
$\psi (x,y)$, using iterated logarithms of $x$ and $y$, is presented.
We establish that
$$
\ln \left (\frac {\psi }{x}\right )=-u
[\ln ^{(2)}x-\ln ^{(2)}y+\ln ^{(3)}x-\ln ^{(3)}y+\ln ^{(4)}x-a] \eqno (2)
$$
where $\underbar{a}<a<\bar{a}$ for some constants
$\underbar{a}$ and $\bar{a}$ (denoting by $\ln ^{(k)}x=\ln ...\ln x$ the $k$-fold iterated
logarithm). The approximation (2) holds in a domain
which is complementary to the one on which the approximation (1) is known to be valid.
One consequence of (2) is an asymptotic expression for Dickman's function, which is of the
form $\ln \rho (u)=-u[\ln u+\ln ^{(2)}u](1+o(1))$, improving known asymptotic
approximations of this type. We employ (2) to establish a version of Bertrand's Conjecture,
indicating how this method may be used to sharpen the result.

\normalsize
\newtheorem{sub}{\name}[section]
\newtheorem{subn}{\name}
\newtheorem{Thm}{Theorem}[section]
\newtheorem{Lem}[Thm]{Lemma}
\newtheorem{Prop}[Thm]{Proposition}
\newtheorem{Cor}[Thm]{Corollary}
\newtheorem{Rem}[Thm]{Remark}
\newtheorem{Def}[Thm]{Definition}
\newtheorem{Ex}[Thm]{Example}
\renewcommand{\thesubn}{}
\newcommand{\dn}[1]{\def\name{#1}}
\newcommand{\be}{\begin{equation}}
\newcommand{\ee}{\end{equation}}
\newcommand{\bs}{\begin{sub}}
\newcommand{\es}{\end{sub}}
\newcommand{\bsn}{\begin{subn}}
\newcommand{\esn}{\end{subn}}
\newcommand{\bea}{\begin{eqnarray}}
\newcommand{\eea}{\end{eqnarray}}
\newcommand{\BA}[1]{\begin{array}{#1}}
\newcommand{\EA}{\end{array}}

\newcommand{\Real}{\mbox{${\rm I\!R}$}}
\newcommand{\real}{{\rm I\!R}}
\newcommand{\Nat}{\mbox{${\rm I\!N}$}}
\newcommand{\thkl}{\rule[-.5mm]{.3mm}{3mm}}
\newcommand{\Proof}{\mbox{\noindent {\bf Proof} \hspace{2mm}}}
\newcommand{\mbinom}[2]{\left (\!\!{\renewcommand{\arraystretch}{0.5}
 \mbox{$\begin{array}[c]{c}
 #1\\ #2
 \end{array}$}}\!\! \right )}
\newcommand{\brang}[1]{\langle #1 \rangle}
\newcommand{\vstrut}[1]{\rule{0mm}{#1mm}}
\newcommand{\rec}[1]{\frac{1}{#1}}
\newcommand{\set}[1]{\{#1\}}
\newcommand{\dist}[2]{\mbox{\rm dist}\,(#1,#2)}
\newcommand{\opname}[1]{\mbox{\rm #1}\,}
\newcommand{\supp}{\opname{supp}}
\newcommand{\mb}[1]{\;\mbox{ #1 }\;}
\newcommand{\undersym}[2]
 {{\renewcommand{\arraystretch}{0.5}
 \mbox{$\begin{array}[t]{c}
 #1\\ #2
 \end{array}$}}}

\newlength{\wex}  \newlength{\hex}
\newcommand{\understack}[3]{%
 \settowidth{\wex}{\mbox{$#3$}} \settoheight{\hex}{\mbox{$#1$}}
 \hspace{\wex}
 \raisebox{-1.2\hex}{\makebox[-\wex][c]{$#2$}}/
 \makebox[\wex][c]{$#1$}
  }%

\newcommand{\smit}[1]{\mbox{\small \it #1}}
\newcommand{\lgit}[1]{\mbox{\large \it #1}}
\newcommand{\scts}[1]{\scriptstyle #1}
\newcommand{\scss}[1]{\scriptscriptstyle #1}
\newcommand{\txts}[1]{\textstyle #1}
\newcommand{\dsps}[1]{\displaystyle #1}

\def\ga{\alpha}     \def\gb{\beta}       \def\gg{\gamma}
\def\gc{\chi}       \def\gd{\delta}      \def\ge{\epsilon}
\def\gth{\theta}                         \def\vge{\varepsilon}
\def\gf{\phi}       \def\vgf{\varphi}    \def\gh{\eta}
\def\gi{\iota}      \def\gk{\kappa}      \def\gl{\lambda}
\def\gm{\mu}        \def\gn{\nu}         \def\gp{\pi}
\def\vgp{\varpi}    \def\gr{\rho}        \def\vgr{\varrho}
\def\gs{\sigma}     \def\vgs{\varsigma}  \def\gt{\tau}
\def\gu{\upsilon}   \def\gv{\vartheta}   \def\gw{\omega}
\def\gx{\xi}        \def\gy{\psi}        \def\gz{\zeta}
\def\Gg{\Gamma}     \def\Gd{\Delta}      \def\Gf{\Phi}
\def\Gth{\Theta}
\def\Gl{\Lambda}    \def\Gs{\Sigma}      \def\Gp{\Pi}
\def\Gw{\Omega}     \def\Gx{\Xi}         \def\Gy{\Psi}
\mysection{Introduction}\label{section1}
A point ${\bf z}$ in $R^m$ is {\em a lattice point} if ${\bf z}=(z_1,...,z_m)$ where each
$z_j$ is an integer. Consider the number of
lattice points included in the simplex $S(a_1,...,a_m)$, where
$$
S(a_1,...,a_m)=\left \{{\bf z}:\sum _{j=1}^m\frac {z_j}{a_j}\leq 1,\,z_j\geq 0,
1\leq j\leq m \right \},
$$
and $a_j$, $j=1,2,...,m$, are positive real numbers. Denote this number by
$\lambda (a_1,...,a_m)$, or $\lambda (S)$.

We need estimates of $\lambda (S)$ as a tool in studying the following problem:
{\em
Let $x$ and $y$ be two positive real numbers, and we are interested in the number of integers
$2\leq k\leq x$ such that the largest prime divisor of $k$ does not exceed $y$, denoted
$\psi(x,y)$.}

Denote by $\{p_j\}_{j=1}^{\infty }$ the
increasing sequence of the primes, and let $m$ be such that
$$
p_{m}<y\leq p_{m+1}.
$$
Then by the Prime Numbers Theorem
\be \label{maprox}
m\approx \frac {y}{\ln y}
\ee
in the sense that the ratio between the two sides of (\ref{maprox}) tends to 1 as
$y\to\infty $. We are thus interested in the integers $k\leq x$ which are of the form
\be \label{kform}
k=\prod _{j=1}^mp_j^{t_j},\,t_j \mbox { are nonnegative integers}.
\ee
Equivalently, we are interested in integers $k$ as in
(\ref{kform}) for which
\be \label{simformula}
\sum _{j=1}^m(\ln p_j)t_j\leq \ln x
\ee
holds. Thus to approximate $\psi (x,y)$ we estimate the expression
$$
\lambda \left (\frac {\ln x}{\ln p_1},...,\frac{\ln x}{\ln p_m}\right ).
$$

There has appeared quite extensive literature on the subject of integers without large
prime divisors since the 30' of the previous century. See e.g. Dickman \cite{Dick},
Erd\H{o}sh \cite{Erd1, Erd2, Erd3}, Erd\H{o}sh and Schinzel \cite{ErdS},
Fouvry and Tenenbaum \cite{Fouv},
Friedlander \cite{Frid1, Frid2, Frid3, Frid4}, Granville \cite{Gran1, Gran2, Gran3, Gran4},
Hazlewood \cite{Hazle}, Hildebrand \cite{Hild1, Hild2, Hild3, Hild4, Hild5},
Hildebrand and Tenenbaum \cite{HildTen1, HildTen2}, Pomerance \cite{Pomer1, Pomer2},
Ramachandra \cite{Ramac1, Ramac2, Ramac3}, Rankin \cite{Rank},  Tenenbaum \cite{Tenen1,Tenen2},
Vershik \cite{Versh}, Xuan \cite{Xu1, Xu2}, and the survey paper by
Hildebrand and Tenenbaum \cite{HildTen3}. More recent related work is presented in
de la Bret\`{e}che and Tenenbaum \cite{Bret},
Hunter \cite{Hunt}, Scourfield \cite{Scour}, Song \cite{Song}, Suzuki \cite{Suzuk1,Suzuk2}
and Tenenbaum \cite{Tenen3},

Dickman \cite{Dick} has established that for every fixed $u>1$ the limit
\be \label{Dickman}
\lim _{x\to \infty }\frac {\psi (x,x^{1/u})}{x}=\rho (u)
\ee
exists, where $\rho (u)$ is  the unique continuous solution of
$$
u\rho '(u)=-\rho (u-1), u\geq 1
$$
satisfying
$$
\rho (u)=1 \mbox { for } 0\leq u\leq 1.
$$
It turns out that $\rho $ satisfies the asymptotic relation
\be \label{asymptrho}
\ln \rho (u)=-(1+o(1))u\ln u.
\ee

Concerning $\psi (x,y)$ we obtained the following result, which is implied by our
main results, Theorems \ref{noaymptot8} and \ref{discsummry}. It deals
with situations where
$$
\ln x<<y<<x,
$$
in a sense expressed precisely Theorem \ref{mainres0}.
We employ the notation
\be \label{iterlog}
\ln ^{(k)}x=\ln \cdots \ln x
\ee
for the $k$th iterated logarithm, where the logarithm function appears $k$ times in the
right hand side of (\ref{iterlog}) and $x$ is sufficiently large. Namely,
$$
\ln ^{(1)}x=\ln x,\,\ln ^{(k+1)}x=\ln (\ln ^{(k)}x),\,k\geq 1.
$$

\begin{Thm} \label{mainres0}
$(i)$ Consider pairs $(x,y)$ such that
\be \label{lwrcondi1}
\exp {(\ln y)^{1-\theta }}<\ln x<\sqrt y
\ee
for some $0<\theta <1$. Denoting
\be \label{unotatn}
u=\frac {\ln x}{\ln y},
\ee
there exist constants $\underline{a}>0$ and $y_0>1$ such that
\be \label{lwrboundassy11n}
\ln \left (\frac {\psi (x,y)}{x}\right )>
-u[\ln ^{(2)}x-\ln ^{(2)}y+\ln ^{(3)}x-\ln ^{(3)}y+\ln ^{(4)}x-\underline{a}]
\ee
for every $y>y_0$.
\newline
$(ii)$
Consider pairs $(x,y)$ such that
\be \label{lwrcondi22}
(\ln y)^{\nu }<\ln x<y^{\beta }
\ee
for some $\nu >2$ and some $0<\beta <1/2$. Then there exist constants
$\overline{a}>a^{\star }$ and $y_0>1$ such that
\be \label{nineqult31n}
\ln \left (\frac {\psi (x,y)}{x}\right )<
-u[\ln ^{(2)}x-\ln ^{(2)}y+\ln ^{(3)}x-\ln ^{(3)}y+\ln ^{(4)}x-\overline{a}]
\ee
for every $y>y_0$.
\end{Thm}

We use the estimates of the iterated logarithms of $x$ and $y$ described in
Proposition \ref{iterlogest} and the inequalities (\ref{lwrboundassy11n}) and
(\ref{nineqult31n}) to obtain the following strengthening of (\ref{asymptrho}).
\begin{Cor}\label{lnroest}
Consider pairs $(x,y)$ such that $(\ref{lwrcondi1})$ holds, and let $u$ be as in
$(\ref{unotatn})$. Then
\be \label{rhouestim}
\ln \rho (u)=-u(\ln u+\ln ^{(2)}u)(1+o(1)),
\ee
where the term $o(1)$ is of order
${\displaystyle O\left (\frac {\ln ^{(3)}u}{\ln ^{(2)}u}\right )}.$
\end{Cor}
Assuming validity of the conjectured expressions (\ref{acumineqlt81}) and (\ref{acumineqlt82})
in Remark \ref{conject} yields that the $o(1)$ term in (\ref{rhouestim}) is of order
${\displaystyle O\left (\frac {\ln ^{(k)}u}{\ln ^{(2)}u}\right )}$ for any $k\geq 3$.

Another application of (\ref{lwrboundassy11n}) and (\ref{nineqult31n}) is to Bertrand's
Conjecture, expressed in Corollary \ref{Bertran}, establishing that for every $\gamma >3/2$
there exists $y_0$ such that
$$
y<p<\gamma y
$$
for some prime $p$, if $y>y_0$. There exist stronger results concerning Bertrand's Conjecture
(see e.g. \cite{Huxl}), and we present Corollary \ref{Bertran} to demonstrate the efficiency
of our main results Theorems \ref{noaymptot8} and \ref{discsummry} as a tool in studying
certain interesting problems.

A uniform version of Dickman's result (\ref{Dickman}) was established by
de Bruijn \cite{Bruij}. Using $u$ in (\ref{unotatn}) he has proved that
\be \label{debruijn}
\psi (x,y)=x\rho (u)\left \{1+O\left (\frac {\ln u}{\ln y}\right )\right \}
\ee
holds uniformly in the domain
$$
2\leq u\leq (\ln y)^{3/5-\epsilon},\,y\geq 2.
$$
This asymptotic relation was extended by Hildebrand \cite{Hild3} who proved that (\ref{debruijn})
holds uniformly in the domain
\be \label{domain2}
2\leq u\leq \exp \{(\ln y)^{3/5-\epsilon}\},\,y\geq 2.
\ee

The upper limit of the domain of validity of (\ref{debruijn}) is related to the error term
in the Prime Number Theorem. Actually Hildebrand established in \cite{Hild1} that
Riemann Hypothesis is true if and only if (\ref{debruijn}) holds uniformly in the domain
\be \label{domain3}
2\leq u\leq y^{1/2-\epsilon},\,y\geq 2,
\ee
for any fixed $\epsilon >0$.

\begin{Rem}\label{compdomain}
Note that restricting to the domain $(\ref{domain3})$, the domain
$$
u<\exp (\ln y)^{3/5-\epsilon }
$$
in $(\ref{domain2})$ is complementary to the domain $$\ln x>\exp {(\ln y)^{1-\theta }}$$
in $(\ref{lwrcondi1})$ for $\theta >2/5+\epsilon $.
\end{Rem}

\begin{Rem}\label{accuracy}
The expressions $(\ref{lwrboundassy11n})$ and $(\ref{nineqult31n})$ provide approximations of
$\ln (\psi (x,y)/x)$, whose accuracy is expressed by $(\overline a-\underline a)u$. To
attain the same level of accuracy as in the approximations $(\ref{debruijn})$ it is required
that $\overline a-\underline a=O(\ln u/\ln x)$, or equivalently
\be \label{adifforder}
\overline a-\underline a=O\left (\frac {\ln \ln x}{\ln x}\right ).
\ee
Moreover, it is easy to see from the proof of Corollary $\ref{Bertran}$, that
$(\ref{adifforder})$ implies the following result for Bertrand's Problem: For every
$\epsilon >0$ there exists a $y_0$ such that if $y>y_0$ then
$$
y<p<y+y^{1/2+\epsilon } \mbox { for some prime } p.
$$
\end{Rem}

We conclude with a result that covers the following range of $(x,y)$
\be \label{lnxlny12}
\frac {1}{2}\ln x<\ln y<\ln x,
\ee
which is different from the ranges indicated in Theorem \ref{mainres0}.
\begin{Thm}\label{niceresult}
Consider the set $E$ of integers $1\leq k\leq x$ for which all the
prime divisors are smaller than $\sqrt x$. $($In our notations $\#(E)=\psi (x,\sqrt x)$.$)$
If $(x,y)$ satisfies $(\ref{lnxlny12})$ then
\be \label{FsqrtN}
\psi (x,y)\geq \#(E)>\alpha x \mbox { for some constant }\alpha >0 \mbox { and every }x>1.
\ee
Actually, for sufficiently large $x$ we may take $\alpha =\ln (e/2)$ in $(\ref{FsqrtN})$.
\end{Thm}
The proof is relegated to the appendix.

The paper is organized as follows. In the next section we
describe a convenient setting for the study of lower and upper bounds of
$\psi(x,y)$. In section \ref{section4} we introduce a family
of auxiliary problems in which our problem can be imbedded. In section \ref{section5}
we introduce our iterations method, which is the main technical tool developed in
this paper. In sections \ref{section6} and \ref{section7} we establish lower and upper bounds
for the auxiliary problems defined in section \ref{section4}, respectively.
Our main results are presented
in section \ref{section8}. In the appendix we establish Theorem \ref{niceresult} and
Proposition \ref{upperbundln2}.

\mysection{The reduced order simplex}\label{section3}
In this section we relate with the high dimensional simplex (\ref{simformula}) a
simplex of smaller order. We will study certain properties of this simplex,
which will be used in the next sections as tools used to establish tight lower
and upper bounds for the number of solutions of (\ref{simformula}).

We divide the integers interval $(1,y)$ into subintervals
\be \label{Jidfn}
J_i=\left (\frac {y}{e^i},\frac {y}{e^{i-1}}\right ), i=1,2,...,r,
\ee
where
\be \label{r1}
r= [\ln y] \mbox { if } \ln y<[\ln y]+\ln 2
\ee
and
\be \label{r2}
r= [\ln y]+1 \mbox { if } \ln y>[\ln y]+\ln 2.
\ee
For simplicity of notations we henceforth consider only case (\ref{r1}),
and comment that the discussion and main results in case (\ref{r2}) are the same.
(In Remark \ref{rdiffrnc} we will indicate where the difference between (\ref{r1})
and (\ref{r2}) plays a role.)

We have for primes $p_j\in J_i$ the relations
\be \label{r3}
\ln y-i<\ln p_j<\ln y-i+1,
\ee
and regarding (\ref{simformula}) this implies
\be \label{Jisum}
(\ln y-i)z_i<\sum _{p_j\in J_i}(\ln p_j)t_j<(\ln y-i+1)z_i,
\ee
where we denote
\be \label{Aidefn}
z_i=\sum _{p_j\in J_i}t_j, \, 1\leq  i\leq r.
\ee
Clearly $(z_1,...,z_r)$ is a nonnegative lattice point in $R^r$.
\begin{Rem}\label{rdiffrnc}
The cases $(\ref{r1})$ and $(\ref{r2})$ differ only when considering $i=r$
in the left inequality of $(\ref{r3})$.
\end{Rem}

If $\{t_j\}_{j=1}^{m}$ is a solution of (\ref{simformula}), then in
view of (\ref{Jisum}) this implies
\be \label{Aiconsrnt}
\sum _{i=1}^r(\ln y-i)z_i<\ln x.
\ee
Therefore the number of solutions
$\{t_j\}_{j=1}^m$ of (\ref{simformula}) is smaller than the number of solutions
$\{t_j\}_{j=1}^m$ of (\ref{Aiconsrnt}).
(We say that $\{t_j\}_{j=1}^m$ is a solution of (\ref{Aiconsrnt}) if
(\ref{Aidefn}) and $(\ref{Aiconsrnt})$ are satisfied.)
Similarly, if $\{t_j\}_{j=1}^m$ is a solution of
\be \label{Aircnt}
\sum _{i=1}^r(\ln y-i+1)z_i<\ln x,
\ee
then in view of (\ref{Jisum}) it is also a solution of (\ref{simformula}),
implying that the number of solutions
$\{t_j\}_{j=1}^m$ of (\ref{simformula}) is larger than the number of solutions
$\{t_j\}_{j=1}^m$ of (\ref{Aircnt}). These considerations are the basis of our
computation of upper and lower bounds for $\psi (x,y)$.

For a prescribed lattice point $(z_1,...,z_r)$ which satisfies (\ref{Aiconsrnt}) we
are interested in the number of lattice points $\{t_j\}_{j=1}^m$ in $R^m$ for which
(\ref{Aidefn}) holds for every $i=1,2,...,r$. Let $m_i$ denote the size of the set
$\{j:p_j\in J_i\}$:
$$
m_i=\#\left \{p_j\in \left (\frac {y}{e^{i}},\frac {y}{e^{i-1}}\right )\right \},
$$
and if $m_i>>1$, then by the Prime Numbers Theorem
\be \label{miapproxi}
m_i\approx \frac {(e-1)y}{(\ln y-i)e^{i}},
\ee
and we have the inequality
\be \label{midefn}
m_i> \frac {y}{e^{i}(\ln y-i)}.
\ee
We denote by $f(k,m)$ the number of different ways in which $k$ can be
written as a sum of $m$ nonnegative integers, and clearly
\be \label{fkdefn}
f(k,m)=\left (\begin{array}{cc}k+m-1\\k\end{array}\right )=
\frac {m(m+1)\cdots (m+k-1)}{k!}.
\ee
Then the number of lattice points $\{t_j\}_{j=1}^m$ that satisfy (\ref{Aidefn})
for every $1\leq i\leq r$ is
\be \label{Kexpress}
K(z_1,...z_r)=\prod _{i=1}^rf(z_i,m_i).
\ee

We Denote by $\overline {\psi }(x,y)$ and $\underline {\psi }(x,y)$ the
number of solutions of (\ref{Aiconsrnt}) and (\ref{Aircnt}) respectively,
and it follows that $\psi (x,y)$ is bounded from above
by  $\overline {\psi }(x,y)$ and from below by $\underline {\psi }(x,y)$.
Using the expression $K(z_1,...,z_r)$ in (\ref{Kexpress}) we consider sums of the form
\be \label{Mexpress}
M(F)=\sum _{{\bf z}\in F}K(z_1,...,z_r),
\ee
where the summation runs over all the lattice points ${\bf z}=\{z_1,...,z_r\}$ which belong to
some set $F$ in $R^r$. Thus when $F$ in (\ref{Mexpress}) is the set of points belonging to
the simplex (\ref{Aiconsrnt}), denoted $F_1$, then by (\ref{fkdefn}) and (\ref{Kexpress})
we have
\be \label{upperbnd}
\overline {\psi}(x,y)=\sum _{\{z_i\}\in F_1}\prod _{i=1}^{r}\frac {m_i^{z_i}}{z_i!}
\left (1+\frac {1}{m_i}\right )\cdots \left (1+\frac {z_i-1}{m_i}\right ).
\ee
Similarly we obtain the following lower bound for $\psi$
\be \label{lowerbnd}
\underline {\psi }(x,y)=\sum _{\{z_i\}\in F_2}\prod _{i=1}^{r}\frac {m_i^{z_i}}{z_i!},
\ee
where $F_2$ is the set of all the lattice points in the simplex (\ref{Aircnt}).

We next consider the product
$$
P_i=\prod _{k=1}^{z_i-1}\left (1+\frac {k}{m_i}\right )
$$
that appears in the right hand side of (\ref{upperbnd}),
and in view of the inequality $\ln (1+t)<t$ for $t>0$ we obtain
$\ln P_i<z_i^2/2m_i$,
hence
\be \label{Pibound}
P_i<e^{z_i^2/2m_i}.
\ee
When dealing with a lower bound we will ignore
the term $\prod _{i=1}^{r}P_i$ in the right hand side of (\ref{upperbnd}), and
we will focus on computing a lower bound to expressions of the form
\be \label{Zexprs}
Z(F)=\sum _{\{z_i\}\in F}\prod _{i=1}^{r}\frac {m_i^{z_i}}{z_i!}
\ee
for certain sets $F$. We will then describe the modifications required to obtain
an upper bound by taking into consideration the terms $P_i$ in (\ref{upperbnd}).

\mysection{A family of auxiliary problems}\label{section4}
It will be convenient to study our main problem, of estimating sums of the form
(\ref{Mexpress}), by using slightly different notations.
In this section we define a collection of problems, parameterized by two real variables,
such that for certain values of the parameters the auxiliary problem coincides with the
main problem. Thus for a positive number $c>1$, let $r=[c]$ and consider the inequality
\be \label{cconstrnt}
cz_0+(c-1)z_1+(c-2)z_2+\cdots +(c-r+1)z_{r-1}<M
\ee
for some positive number $M>1$, where ${\bf z}=\{z_i\}_{i=0}^{r-1}$ is a nonnegative
lattice point in  $R^{r}$ (compare with (\ref{Aircnt})).
We associate with $c$ the $r$ bases
\be \label{mbases}
m_i=\frac {(e-1)e^{c-i}}{c-i},\,0\leq i\leq r-1
\ee
(compare with  (\ref{miapproxi}) in case that $c=\ln n$). In view of (\ref{lowerbnd})
we address the problem of computing the sum
\be \label{cMsum}
F(c,M)=\sum _{\bf z}\prod _{i=0}^{r-1}\frac {m_i^{z_i}}{z_i!},
\ee
where ${\bf z}=(z_0,...,z_{r-1})$ runs over all the nonnegative lattice points which satisfy
(\ref{cconstrnt}); we call this {\em Problem $P_{c,M}$} for the
$r$ variables $z_0$,...,$z_{r-1}$.
\begin{Rem}\label{ourapplic}
There is a close relation between the value of Problem $P_{c,M}$ and $\psi (x,y)$ for
\be \label{ourvalus}
c=\ln y\mbox { and } M=\ln x.
\ee
Thus the value of $P_{c,M}$ yields a lower bound for $\psi(x,y)$. We also note that
if $c\geq M$ $($namely $y\geq x)$ and $x$ is an integer, then
\be \label{nuNMeql}
\psi (x,y)=x=e^M.
\ee
\end{Rem}

To establish an upper bound for $\psi (x,y)$ we will estimate a sum of the type
(\ref{Mexpress}), which is associated with the simplex
\be \label{inter}
(c-1)z_1+(c-2)z_2+\cdots +(c-r)z_{r}<M
\ee
(compare with (\ref{Aiconsrnt})). This sum is smaller than the corresponding sum that
is associated with the simplex
\be \label{ccconstrnt}
cz_0+(c-1)z_1+(c-2)z_2+\cdots +(c-r)z_{r}<M,
\ee
which we denote by $G_0(c,M)$. Thus to obtain an upper bound for $G_0(c,M)$
we consider a sum similar to the one in (\ref{cMsum}), where
we take into consideration the
terms $P_i$ in (\ref{Pibound}). We then address the problem of computing the sum
\be \label{GcMsum}
G(c,M)=\sum _{\bf z}\prod _{i=0}^{r}\frac {m_i^{z_i}e^{z_i^2/m_i}}{z_i!},
\ee
where ${\bf z}=(z_0,z_1,...,z_{r})$ runs over all the nonnegative lattice points which
satisfy (\ref{ccconstrnt}); we call this {\em Problem $Q_{c,M}$} for the
$r+1$ variables $z_0$,$z_1$,...,$z_{r}$.
\begin{Rem}
We use the simplex $(\ref{ccconstrnt})$ rather than the simplex $(\ref{inter})$, which is
more directly related to $(\ref{Aiconsrnt})$, in order to avoid repetition of
computations for the lower and upper bounds. Thus a substantial part of the computations
for $(\ref{cconstrnt})$ and $(\ref{ccconstrnt})$ will be unified.
\end{Rem}

We claim that for a fixed value of $z_0$, Problem $P_{c,M}$ reduces to
Problem $P_{c-1,M-cz_0}$ for the $r-1$ variables $z_1$,...,$z_{r-1}$. To justify
this statement we have to check that the $r-1$ bases $m_1$,...,$m_{r-1}$ in
(\ref{mbases}) are indeed the bases associated with Problem $P_{c-1,M-cz_0}$,
which is easily verified.

The possible values for the variable $z_0$ in (\ref{cconstrnt}) are the
integers $z$ satisfying
$$
0\leq z\leq \frac {M}{c},
$$
and it follows from (\ref{cMsum}) that
\be \label{Frecurs}
F(c,M)=\sum _{z=0}^{[M/c]}F(c-1,M-cz)\frac {m_0^z}{z!}.
\ee

In the subsequent discussion we will consider situations where $F(\cdot ,\cdot )$
satisfies inequalities of the form
\be \label{Fform}
F(c,M)\geq Be^{M\left (1-\frac {\ln M}{c+1}
+\frac {\gamma }{c+1}\right )}
\ee
for some constant  $0<B\leq 1$. In terms of the original parameters we are actually
interested in inequalities of the form
\be \label{nuform}
\psi (x,y)\geq Bx^{\left (1-\frac {\ln \ln x}{\ln y+1}
+\frac {\gamma }{\ln y+1}\right )},
\ee
where $(x,y)$ and $(c,M)$ are related as in (\ref{ourvalus}).
\begin{Rem} \label{Mc2}
We note that $M/c$ is the parameter $u$ in $(\ref{unotatn})$, which appeared, e.g., in
$(\ref{debruijn})$, $(\ref{domain2})$ and $(\ref{domain3})$. It follows from $(\ref{FsqrtN})$ 
in Theorem $\ref{niceresult}$ that for a fixed $\gamma $, inequality $(\ref{nuform})$
holds whenever $M/c<2$. Indeed, for $M=\ln x$ and
$c=\ln y$ the condition $M/c<2$ translates to $y>\sqrt x$, implying
$\psi (x,y)>\alpha x$ by $(\ref{FsqrtN})$. But the inequality
$$
\alpha x>x^{1-\frac {\ln \ln x}{\ln y+1}+\frac {\gamma }{\ln y+1}}
$$
is equivalent to
$$
\frac {\ln x}{\ln y+1}(\ln \ln x-\gamma )>-\ln \alpha ,
$$
and this holds for every $x>x_0$, for some $x_0$, since $y<x$. For $x\leq x_0$,
however, $(\ref{nuform})$ holds for some $B(\gamma )$,
since in this case we have a bounded set of pairs $(x,y)$.
Therefore, when trying to establish an inequality of the type $(\ref{Fform})$,
we may assume that
\be \label{Mclrg2}
\frac {M}{c}\geq 2,
\ee
since for $M/c<2$ inequality $(\ref{nuform})$ is already established.
\end{Rem}
\mysection{The iterations method}\label{section5}
The discussion in this section is fundamental to our analysis. We develop the iterations
method which will be employed in the subsequent sections to establish lower and upper
bounds for $\psi $.

Assume that for a certain $\gamma >0$ and some $0<B<1$,
inequality (\ref{Fform}) holds for any pair $(c,M)$ which verifies
\be \label{cMM0cond}
c\leq \kappa _0,
\ee
for a certain $\kappa _0$. We consider then pairs $(c,M)$ that satisfy
\be \label{nextblock1}
\kappa _0<c\leq \kappa _0+1,
\ee
and our goal is to establish the inequality $(\ref{Fform})$ for such pairs as well. Once
this is achieved we will iterate the argument to obtain a lower bound for all pairs in a
certain domain.

Intending to employ (\ref{Frecurs}) to establish a lower bound to $F(c,M)$, and assuming
that (\ref{Fform}) holds whenever (\ref{cMM0cond}) is satisfied, we
will estimate from below the expressions
\be \label{Fzbarexp}
F(c-1,M-cz)\frac {m_0^{z}}{z!}
\ee
which appear in (\ref{Frecurs}), and this
for integers $0\leq z\leq M/c$. By (\ref{nextblock1}) $c-1\leq \kappa _0$, and
we may use (\ref{Fform}) for the pair $(c-1,M-cz)$, obtaining
\be \label{c1Mzbnd}
F(c-1,M-cz)\geq Be^A,
\ee
where
\be \label{Alogaritm}
A=(M-cz)-\frac {1}{c}(M-cz)\ln (M-cz)+\frac {(M-cz)\gamma}{c}.
\ee
Moreover, the inequality
\be \label{mzeE}
\frac {m_0^z}{z!}>e^E,
\ee
holds, where we denote
\be \label{Elogaritm1}
E=\left (z\ln m_0-z\ln z+z\right )-\left (\frac {1}{2}\ln z
+\frac {1}{2}\ln \pi +\frac {3}{2}\ln 2\right ),
\ee
where we used Stirling's formula
\be \label{Stirapprx}
St(z)=\sqrt {2\pi z}\left (\frac {z}{e}\right )^z
\ee
to estimate
\be \label{2gamaz}
z!<2St(z) \mbox { for every } z\geq 1.
\ee
In (\ref{Elogaritm1}), a term $(-\ln 2)$  arises from the factor $2$ in (\ref{2gamaz}),
and the term
\be \label{dist}
-\frac {1}{2}(\ln z+\ln \pi +\ln 2)
\ee
is due to the logarithm of $\sqrt {2\pi z}$ in (\ref{Stirapprx}).
To avoid the disturbing term (\ref{dist}) in (\ref{Elogaritm1}) we note that
\be \label{zpibeta}
\frac {1}{2}\ln z+\frac {1}{2}\ln \pi +\frac {3}{2}\ln 2<hz
\ee
where $h>0$ may be chosen arbitrarily
small provided that $z$ is sufficiently large. It follows that
\be \label{zineqmodi}
z-\left (\frac {1}{2}\ln z
+\frac {1}{2}\ln \pi +\frac {3}{2}\ln 2\right )>bz
\ee
where
\be \label{b1beta}
b=1-h
\ee
may be chosen arbitrarily close to 1
provided that $z$ is large enough, and we thus obtain
\be \label{Elogaritm}
E>\left (z\ln m_0-z\ln z+bz\right )
\ee
for sufficiently large values of $z$.

It follows from $m_0=(e-1)e^{c}/c$ that
$$
z\ln m_0=cz-z\ln c+z\ln (e-1).
$$
Using the last equation in (\ref{Elogaritm}) and recalling (\ref{Alogaritm}) yield that
\be \label{AEestim}
A+E>H(z),
\ee
denoting
\be \label{maxexpres}
H(z)=M\left (1+\frac {\gamma }{c}\right )+(a-\gamma )z-\frac {M}{c}\ln c
-z\ln z-\left (\frac {M}{c}-z\right )\ln \left (\frac {M}{c}-z\right )
\ee
and
\be \label{adefinitn}
a=b+\ln (e-1).
\ee
Thus $a$ is smaller and arbitrarily close to $a^{\star }$, which is defined by
\be \label{limita}
a^{\star }=1+\ln (e-1).
\ee
It follows from (\ref{c1Mzbnd}), (\ref{mzeE}) and (\ref{AEestim}) that
\be \label{FcMHD}
F(c-1,M-cz)\frac {m_0^z}{z!}>Be^{H(z)},
\ee
and to obtain a lower bound for the sum in (\ref{Frecurs}) we will estimate
the maximal value of $H(z)$, $0\leq z\leq [M/c]$, where $z$ is an integer.

\begin{Rem}\label{newrem}
We will compute a maximizer $z_0$ of $H(\cdot )$ defined on the real interval $[0,[M/c]]$,
and in general $z_0$ is not an integer. Let $z_1$ be the integer
$$
z_1=z_0+\theta \mbox { for some }
0\leq \theta <1,
$$
and then
$$
H(z_1)=H(z_0)+\frac {1}{2}H''(\zeta )\theta ^2
$$
for some $z_0<\zeta <z_1$. But
$$
H''(\zeta )=\frac {-M/c}{\zeta (M/c-\zeta )},
$$
and it follows from $\zeta \geq 1$ that
$$
|H''(\zeta )|\leq \frac {M/c}{M/c-1}<2
$$
$($since $M/c>2)$, and we obtain
\be \label{Hz0z1}
H(z_1)>H(z_0)-\theta ^2.
\ee
Similarly, for the integer $z_2=z_0-(1-\theta )$ we have
\be \label{Hz0z2}
H(z_2)>H(z_0)-(1-\theta )^2.
\ee
It follows from $(\ref{FcMHD})$, $(\ref{Hz0z1})$ and $(\ref{Hz0z2})$ that
\be \label{intgrestm}
\sum _{z=0}^{[M/c]}F(c-1,M-cz)\frac {m_0^z}{z!}>B\left (e^{H(z_1)}+e^{H(z_2)}\right )
>Be^{H(z_0)}
\ee
since
$$
\min _{0\leq \theta \leq 1}\left \{e^{-\theta ^2}+e^{-(1-\theta )^2}\right \}>1.
$$
Therefore we may use the maximal value of $H(z)$ over the whole real interval
$0\leq z\leq M/c$.
\end{Rem}
We have the following basic result.
\begin{Prop}\label{maxHvalue}
Let $H(z)$ be as in $(\ref{maxexpres})$. Then
\be \label{Hzmaxlower1}
\max \left \{H(z):0\leq z\leq \frac {M}{c}\right \}=
M\left (1-\frac {\ln M}{c}+\frac {\gamma +f(\gamma )}{c}\right ),
\ee
where
\be \label{fdefn}
f(\gamma )=\ln (1+e^{a-\gamma }).
\ee
\end{Prop}
{\em Proof}:
Denoting
$$
u=\frac {M}{c} \mbox { and }z=ut
$$
it follows that
\begin{eqnarray}\label{zmaxval}
\max _z\{(a-\gamma )z-z\ln z-(u-z)\ln (u-z)\}=\nonumber \\
-u\ln u+u\max _{0\leq t\leq 1}\{(a-\gamma )t-t\ln t-(1-t)\ln (1-t)\}.
\end{eqnarray}
We denote
\be \label{varphidfn}
\varphi (t)=(a-\gamma )t-t\ln t-(1-t)\ln (1-t),
\ee
and it follows that the maximizer $t_0$ of $\varphi $ satisfies
$$
(a-\gamma )-\ln t_0+\ln (1-t_0)=0.
$$
We conclude that
\be \label{toexpr}
t_0(\gamma )=\frac {1}{1+e^{\gamma -a}},
\ee
and the maximal value of $\varphi (\cdot )$ is given by
$$
(a-\gamma )t_0+\ln (1+e^{\gamma -a})-(1-t_0)(\gamma -a),
$$
which yields
\be \label{tmaxval}
\max \{\varphi (t):0\leq t\leq 1\}=\ln (1+e^{a-\gamma }).
\ee
We thus conclude from (\ref{zmaxval}) and (\ref{tmaxval}) that
\begin{eqnarray} \label{Hzmaximal}
\max _{0\leq z\leq u}\left \{z(a-\gamma )
-z\ln z-\left (\frac {M}{c}-z\right )
\ln \left (\frac {M}{c}-z\right )\right \}
=\nonumber \\ -u\ln u
+u\ln (1+e^{a-\gamma }).
\end{eqnarray}
It follows from (\ref{maxexpres}) and (\ref{Hzmaximal}) that (\ref{Hzmaxlower1}) is
satisfied, where $f(\gamma )$ in (\ref{fdefn}) is the maximum
in (\ref{tmaxval}). The proof of the proposition is complete. $\hfill \Box$
\begin{Prop}\label{formua430}
Assume that
$$
F(c',M)\geq Be^{M\left (1-\frac {\ln M}{c'+1}
+\frac {\gamma }{c'+1}\right )}
$$
holds for every $c'\leq c-1$, for some $c>1$. Then
\be \label{Hzmaxlower}
F(c,M)\geq B\exp \left \{M\left (1-\frac {\ln M}{c}
+\frac {\gamma+f(\gamma )}{c}\right )\right \}.
\ee
\end{Prop}
{\em Proof}: Equation (\ref{Hzmaxlower})
follows from (\ref{Frecurs}), (\ref{intgrestm}) and (\ref{Hzmaxlower1}). $\hfill \Box$

For the induction argument we need that (\ref{Fform}) would hold for some initial value of
$c$, say for $c=\kappa $ for some $\kappa >1$. This is the content of the following result.
\begin{Prop}\label{Bkapagama}
For a prescribed $\gamma >0$ the inequality
\be \label{Fkapgamform}
F(\kappa ,M)\geq B(\kappa ,\gamma )e^{M\left (1-\frac {\ln M}{\kappa +1}
+\frac {\gamma }{\kappa +1}\right )}
\ee
holds for every $M\geq 0$, where
\be \label{Bkapgamval}
B(\kappa ,\gamma )=e^{-e^{\kappa +\gamma }}.
\ee
\end{Prop}
{\em Proof}: The maximal value of
$$
M\mapsto M\left (1-\frac {\ln M}{\kappa +1}+\frac {\gamma }{\kappa +1}\right )
$$
is $ \displaystyle {\frac {e^{\kappa +\gamma }}{\kappa +1}}$, and it is attained at
$M_0=e^{\kappa +\gamma }$. Since $B(\kappa ,\gamma )$ in (\ref{Bkapgamval}) satisfies
$$
B(\kappa ,\gamma )e^{\frac {e^{\kappa +\gamma }}{\kappa +1}}<1,
$$
and since $F(c,M)\geq 1$, inequality (\ref{Fkapgamform}) follows for every $M>1$.
$\hfill \Box$.

It follows from Proposition \ref{Bkapagama} that if $B$ in (\ref{Fform}) is equal to
$B(\kappa ,\gamma )$  in (\ref{Bkapgamval}), then (\ref{Fform}) holds for any pair
$(c,M)$ such that $c\leq \kappa$.
\mysection{A lower bound for Problem $P_{c,M}$}\label{section6}
In this section we employ the results of the previous section to establish a lower bound
for Problem $P_{c,M}$. We will construct a sequence
$$
\{(c_j,M_j)\}_{j=0}^l
$$
(where $c_{j+1}=c_j-1$), for which
(\ref{Hzmaxlower}) will be employed successively. The coefficient $B$ will be chosen
such that
\be \label{mainesti}
F(c,M)\geq B\exp \left \{M\left (1-\frac {\ln M}{c+1}
+\frac {\gamma '}{c+1}\right )\right \}
\ee
will hold for the pair $(c_l,M_l)$ for a certain $\gamma '=\gamma _l$, and consequently,
employing (\ref{Hzmaxlower}),
it will hold for each $(c_j,M_j)$  with a corresponding $\gamma '=\gamma _j$. In
particular it will hold for $(c,M)=(c_0,M_0)$ with a certain $\gamma '=\gamma _0$.

Recall that in deriving the estimate (\ref{Hzmaxlower}) we used a value
$$
z_0=ut_0
$$
and that we associated with $(c_0,M_0)$ a pair $(c_1,M_1)$ such that $c_1=c_0-1$, and
\be \label{M1expr}
M_1=M_0(1-t_0).
\ee
Although this pair does not correspond to an integer $z$, it may be used in computing
a lower bound for $F(c,M)$, as explained in Remark \ref{newrem}.

Let $a$ be associated with $z_0$ as in $(\ref{zpibeta})$, $(\ref{b1beta})$ and
$(\ref{adefinitn})$. Recalling (\ref{limita}) we have the following result:
\begin{Prop}\label{aureltn}
For any prescribed $\epsilon >0$ there exists a $u_0$ such that
\be \label{astaraclose}
|a-a^{\star }|<\epsilon \mbox { if } u>u_0.
\ee
\end{Prop}
Concerning (\ref{Hzmaxlower}), we wish to estimate its right hand side as follows:
\be \label{Mgamampr}
M\left (1-\frac {\ln M}{c}
+\frac {\gamma+f(\gamma )}{c}\right )>
M\left (1-\frac {\ln M}{c+1}
+\frac {\gamma '}{c+1}\right )
\ee
for a certain $\gamma '$. Clearly the inequality (\ref{Mgamampr}) is equivalent to
\be \label{Mgamampr1}
\frac {\gamma +f(\gamma )}{c}>\frac {\ln M}{c(c+1)}+\frac {\gamma '}{c+1}.
\ee
For any $\beta >0$ we denote
\be \label{Dbetadfn}
{\cal D}_{\beta }=\{(c,M):1\leq M\leq e^{\beta c}\},
\ee
and for a fixed $\alpha >0$ and a pair $(c,M)$ we denote
\be \label{gamcM}
\gamma _{c,M}=a-\alpha +\ln c-\ln \ln M.
\ee

We assume the validity of (\ref{mainesti}) with $c-1$ replacing $c$ and with
$$
\gamma '=\gamma _{c-1,M-cz}.
$$
Namely we assume that
\be \label{Fc1lowbnd}
F(c-1,M')\geq B \exp \left \{M'\left (1-\frac {\ln M'}{c}
+\frac {\gamma _{c-1,M'}}{c}\right )\right \}
\ee
for every $1\leq M'\leq M$. Using (\ref{gamcM}) in (\ref{Fc1lowbnd}) yields
$$
F(c-1,M')\geq B \exp \left \{M'\left (1-\frac {\ln M'}{c}
+\frac {a-\alpha +\ln (c-1)-\ln \ln M}{c}\right )\right \},
$$
which we write in the form
\be \label{Fc1M'gam0}
F(c-1,M')\geq B \exp \left \{M'\left (1-\frac {\ln M'}{c}
+\frac {\gamma _0}{c}\right )\right \}
\ee
for every $1\leq M'\leq M$, denoting
\be \label{gama0expr}
\gamma _0=a-\alpha +\ln (c-1)-\ln \ln M.
\ee

For a pair $(c,M)$ we consider the maximization over $z$ of
\be \label{maxc1Fczm0}
F(c-1,M-cz)\frac {m_0^z}{z!}.
\ee
The fact that the parameter $\gamma _0$ in (\ref{Fc1M'gam0}) is one and the same for all $M'$
enables to employ the results of section \ref{section5}. Thus the maximal value of
(\ref{maxc1Fczm0}) exceeds the maximal value which is obtained when we replace $F(c-1,M-cz)$
by the right hand side of (\ref{Fc1M'gam0}), with $M'=M-cz$, namely the maximal value of
\be \label{boundsmaxi}
B\exp \left \{(M-cz)\left [1-\frac {\ln (M-cz)}{c}
+\frac {\gamma _0}{c}\right ]\right \}\frac {m_0^z}{z!}
\ee
over $0\leq z\leq M/c$. This latter maximum is attained at
\be \label{M'Mdfn}
M'=M(1-t_0)
\ee
where
\be \label{t0fdormul}
t_0=\frac {1}{1+e^{\gamma _0-a}}.
\ee
We focus our attention on the domain ${\cal D}_{1/2}$ (recall (\ref{Dbetadfn})), and will next
establish that if $(c,M)\in {\cal D}_{1/2}$ then also the resulting pair $(c-1,M')$ belongs
to ${\cal D}_{1/2}$.
\begin{Prop}\label{Dbetacontn}
There exists an $\alpha _0>0$ and $c_0>0$ such that if $\alpha $ in $(\ref{gamcM})$
satisfies $\alpha >\alpha _0$  then for $c>c_0$
\be \label{betaclose}
(c,M)\in {\cal D}_{1/2}\Rightarrow (c-1,M')\in {\cal D}_{1/2}.
\ee
\end{Prop}
{\em Proof}:  By (\ref{gama0expr})
\be \label{expgam0a}
e^{\gamma _0-a}=e^{-\alpha }\frac {c-1}{\ln M},
\ee
and since $(c,M)\in {\cal D}_{1/2}$, we have $(\ln M )/c\leq 1/2$.
We distinguish between the situation where $(\ln M )/c$ is close to $1/2$, and where
$(\ln M )/c$ is smaller, say
\be \label{lnMucval}
\frac {\ln M}{c}<\mu
\ee
for some $0<\mu <1/2$. If (\ref{lnMucval}) holds then for some $c_0$ we have
$$
\frac {\ln M'}{c-1}<\frac {\ln M}{c-1} <\frac {1}{2}
$$
for every $c>c_0$. If, however, (\ref{lnMucval}) does not hold,
so that
\be \label{lnMcval}
\frac {c}{\ln M}\leq \frac {1}{\mu},
\ee
then we obtain from (\ref{t0fdormul}) and (\ref{expgam0a}) that $t_0<1$ is arbitrarily
close to 1 provided that $\alpha $ is large enough. In particular we have that
$$
-\ln (1-t_0)>\frac {1}{2},
$$
which implies $(\ln M')/(c-1)<1/2$ in view of $\ln M'=\ln M+\ln (1-t_0)$.
The proof is complete.
$\hfill \Box$

\begin{Rem}\label{smallbeta}
We consider pairs $(c,M)\in {\cal D}_{\beta }$ where we let $\beta \to 0$. It then follows
from $(\ref{gama0expr})$ and $(\ref{t0fdormul})$ that
\be \label{frstineq}
\ln (1-t_0)<-ke^{a-\gamma _0}
\ee
for some $0<k<1$. Actually $k$ is arbitrarily close to $1$ if $\beta $ is sufficiently small,
since then, by $(\ref{gama0expr})$, $\gamma _0$ becomes arbitrarily large, using
$$
\ln c-\ln \ln M\approx \ln (1/\beta ).
$$
It follows from $(\ref{gama0expr})$ that
\be \label{secndineq}
e^{a-\gamma _0}=e^{\alpha }\frac {\ln M}{c-1},
\ee
and employing
$$
\frac {\ln M'}{c-1}=\frac {\ln M+\ln (1-t_0)}{c-1}
$$
we conclude from $(\ref{frstineq})$ and $(\ref{secndineq})$ that we have
\be \label{5star}
\frac {\ln M'}{c-1}<\frac {\ln M(c-1-ke^{\alpha })}{c(c-1)}<\frac {\ln M}{c}<\beta
\ee
for $\alpha >\alpha (\beta )$, where $\alpha (\beta )\to 0$ as $\beta \to 0$, and
actually we may take
\be \label{alphabeta}
\alpha (\beta )=\frac {\beta }{2}\left (1+\frac {\beta }{2}\right ).
\ee
Thus for sufficiently small $\beta $ we have the implication
\be \label{betaimplictn}
(c,M)\in {\cal D}_{\beta }\Rightarrow (c-1,M')\in {\cal D}_{\beta }.
\ee
\end{Rem}

We will next establish (\ref{mainesti}) with
\be \label{gamaprchoic}
\gamma '=\gamma _{c,M}
\ee
(recall (\ref{gamcM})), assuming the validity of (\ref{mainesti}) with $c$ being replaced
by $c-1$.
\begin{Prop}\label{iteratineq}
Let $z_0$ be the maximizer in the maximization over $z$ of $(\ref{boundsmaxi})$, and let
$a$ be associated with $z_0$ as in $(\ref{zpibeta})$, $(\ref{b1beta})$ and
$(\ref{adefinitn})$. Let $\gamma '$ be as in $(\ref{gamaprchoic})$ and $\gamma =\gamma _0$
$($recall $(\ref{gama0expr}))$, and assume that $(c,M)\in {\cal D}_{1/2}$.
Then $(\ref{Mgamampr})$ holds for some $\alpha >0$.
\end{Prop}
{\em Proof}: We consider the expression
\be \label{fgamacM}
f(\gamma )=f(\gamma _0)=\ln \left (1+e^{\alpha }\frac {\ln M}{c-1}\right ).
\ee
Let $K>0$ be such that
$$
\ln (1+Kx)>2x \mbox { for every } 0<x<1/2
$$
(e.g. $K=4$). We then have that for some $\alpha _0>0$
\be \label{ln1beta}
\ln \left (1+e^{\alpha }\frac {\ln M}{c-1}\right )>2\frac {\ln M}{c}
\ee
for every $\alpha >\alpha _0$, if $(c,M)\in {\cal D}_{1/2}$. It then follows from
(\ref{fgamacM}) that
$$
f(\gamma )>2\frac {\ln M}{c},
$$
and to establish (\ref{Mgamampr1}) it is enough to verify
\be \label{gmgmprm}
\frac {\gamma }{c}-\frac {\gamma '}{c+1}>-\frac {\ln M}{c}.
\ee
In view of
$$
\gamma =a-\alpha +\ln (c-1)-\ln \ln M
$$
and
$$
\gamma '=a-\alpha +\ln c-\ln \ln M
$$
we have
$$
\frac {\gamma }{c}-\frac {\gamma '}{c+1}=\frac {a-\alpha -\ln \ln M}{c(c+1)}
+\frac {\ln (c-1)}{c}-\frac {\ln c}{c+1},
$$
and (\ref{gmgmprm}) follows from
$$
\frac {\ln (c-1)}{c}-\frac {\ln c}{c+1}>-\frac {1}{c(c-1)}.
$$
This concludes the proof of the proposition. $\hfill \Box$

\begin{Rem}\label{alpbetrel}
Suppose that rather than $(c,M)\in {\cal D}_{1/2}$ we consider
$(c,M)\in {\cal D}_{\beta }$ with $\beta <<1$. Arguing as in the above proof and
analogous to $(\ref{ln1beta})$ we consider the inequality
\be \label{ln1beta1}
\ln \left (1+e^{\alpha }\beta \right )>(1+\epsilon )\beta
\ee
for arbitrarily small $\epsilon >0$. For $\epsilon =0$ let $(\ref{ln1beta1})$ hold
for $\alpha >\alpha (\beta )$, and it is easy to see that we may take $\alpha (\beta )$
as in $(\ref{alphabeta})$.
\end{Rem}

It follows from (\ref{Hzmaxlower}), (\ref{Mgamampr}) and Proposition \ref{iteratineq}
that the inequality
\be \label{lowerboundcM}
F(c,M)>B\exp \left \{M\left (1-\frac {\ln M}{c+1}+\frac {a-\alpha +\ln c-\ln \ln M}{c+1}
\right )\right \}
\ee
holds for certain values of $\alpha $ and certain pairs $(c,M)$. Actually, the above
discussion yields the next iterative property.
\begin{Prop}\label{kapkap1lwr}
There exists an $\alpha _0>0$ such that
for any fixed $\alpha >\alpha _0$ there exists $\kappa _0>0$ with the following property:
If $\kappa >\kappa _0$ is such that  $(\ref{lowerboundcM})$ holds for every
$(c,M)\in {\cal D}_{1/2}$ satisfying $\kappa _0<c\leq \kappa $, then
$(\ref{lowerboundcM})$ also holds for every $(c,M)$ that verifies
$$
(c,M)\in {\cal D}_{1/2} \mbox { and } \kappa _0<c\leq \kappa +1.
$$
\end{Prop}

To start the iterations procedure we need the following result:
\begin{Prop}\label{lwrboundres}
For a fixed $\alpha >\alpha _0$ let $\kappa _0$ be as in Proposition
$\ref{kapkap1lwr}$, and let $B$ be defined by
\be \label{Bkap0}
B=e^{-\kappa _0e^{a+\kappa _0}}.
\ee
Then $(\ref{lowerboundcM})$ holds for every $(c,M)\in {\cal D}_{1/2}$ such that
$c\geq \kappa _0$.
\end{Prop}
{\em Proof}: The assertion of the proposition follows from Propositions \ref{Bkapagama}
and \ref{kapkap1lwr}, employing an induction argument. $\hfill \Box$

We conclude from Propositions \ref{kapkap1lwr} and \ref{lwrboundres} the following result.
\begin{Prop}\label{noaymptot}
Let $a<a^{\star }$ be fixed. Then there exist $c_0$, $\alpha $ and $B$ such that
\be \label{lwrboundassy}
F(c,M)>B\exp \left \{M\left (1-\frac {\ln M+\ln \ln M}{c+1}+\frac {a-\alpha +\ln c}{c+1}
\right )\right \}
\ee
for every $(c,M)$ such that $M<e^{c/2}$ and $c>c_0$.
\end{Prop}

We next consider the expression
\be \label{gamccMM}
\bar \gamma _{c,M}=a-\alpha +\ln c+\ln \ln c-\ln \ln M-\ln \ln \ln M
\ee
instead of the expression $\gamma _{c,M}$ in (\ref{gamcM}), and repeat the above
argument and computation using $\bar \gamma _{c,M}$ rather than $\gamma _{c,M}$. We
will next indicate the required modifications.

Instead of (\ref{gama0expr}) we have now
\be \label{gama0new}
\gamma _0=a-\alpha +\ln (c-1)+\ln ^{(2)}(c-1)-\ln ^{(2)}M-\ln ^{(3)}M.
\ee
Proposition \ref{Dbetacontn} and its proof still hold, where instead of (\ref{expgam0a})
we have now
\be \label{expgam0anew}
e^{\gamma _0-a}=e^{-\alpha }\frac {c-1}{\ln M}\frac {\ln (c-1)}{\ln \ln M}.
\ee
We note that if (\ref{lnMcval}) holds then $\ln c<\ln \ln M+\ln (1/\mu )$, implying
that
$$
\frac {\ln (c-1)}{\ln \ln M}<\frac {3}{2} \mbox { if } c>c_0,
$$
for some $c_0>0$. The rest of the proof of Proposition \ref{Dbetacontn} applies in the
present case without change.

Concerning the proof of Proposition \ref{iteratineq}, using the expression
(\ref{gama0new}) for  $\gamma _0$, we obtain
\be \label{fgam0new}
f(\gamma _0)=\ln \left (1+e^{\alpha }\frac {\ln M}{c-1}\frac {\ln \ln M}{\ln (c-1)}
\right ).
\ee
We note that by $\ln M<c/2$ we have $\ln ^{(2)} M/\ln (c-1)<1$. Moreover, assuming that
\be \label{ctetlnM}
c^{1-\theta }<\ln M
\ee
for some $0<\theta <1$ we obtain
\be \label{lnlnMlncest}
\frac {\ln \ln M}{\ln c}>1-\theta.
\ee
Using (\ref{lnlnMlncest}) in (\ref{fgam0new}) we can employ the rest of the proof of
Proposition \ref{iteratineq} to establish the following result.
\begin{Prop}\label{iteratinew}
Let $z_0$ be the maximizer in the maximization over $z$ of $(\ref{boundsmaxi})$, and let
$a$ be associated with $z_0$ as in $(\ref{zpibeta})$, $(\ref{b1beta})$ and
$(\ref{adefinitn})$. Let
\be \label{gammanew}
\gamma =a-\alpha +\ln (c-1)+\ln \ln (c-1)-\ln \ln M-\ln \ln \ln M
\ee
and
\be \label{gammaprnew}
\gamma '=a-\alpha +\ln c+\ln \ln c-\ln \ln M-\ln \ln \ln M,
\ee
and assume that $(\ref{ctetlnM})$ and $(c,M)\in {\cal D}_{1/2}$ are satisfied.
Then there exists an $\alpha $ such that $(\ref{Mgamampr})$ holds.
\end{Prop}

The following is the lower bound which we obtain for $F(c,M)$.
\begin{Thm} \label{mainlwrbndrslt}
Consider pairs $(c,M)$ such that
\be \label{lwrcondi}
c^{1-\theta}<\ln M<\frac {1}{2}c
\ee
for some $0<\theta <1$, and let $a<a^{\star }$ be fixed. Then there exist constants
$\alpha >0$ and $c_0>1$ such that
\be \label{lwrboundassy1}
F(c,M)>\exp \left \{M\left (1-\frac {\ln M+\ln ^{(2)}M+\ln ^{(3)}M}{c+1}
+\frac {a-\alpha +\ln c+\ln ^{(2)}c}{c+1}
\right )\right \}
\ee
for every $c>c_0$.
\end{Thm}
\begin{Rem}\label{betateta}
In view of Remarks $\ref{smallbeta}$ and $\ref{alpbetrel}$
the discussion and proof which yield Theorem $\ref{mainlwrbndrslt}$ can be employed to conclude
the following: for any $\alpha >0$, which may be arbitrarily small, we can choose $\beta >0$
and $\theta >0$ sufficiently small such that $(\ref{lwrboundassy1})$ holds for pairs $(c,M)$
satisfying
\be \label{modlwrcondi}
c^{1-\theta}<\ln M<\beta c,
\ee
replacing $(\ref{lwrcondi})$. Actually, in view of $(\ref{alphabeta})$, we may take
$\beta <(2-\epsilon )\alpha $, if $\alpha $ is small enough.
\end{Rem}
\mysection{An upper bound for Problem $Q_{c,M}$}\label{section7}
In this section we are concerned with the upper bound for $G(c,M)$ in (\ref{GcMsum}).
We will employ a method similar to the one used to establish a lower bound
for $F(c,M)$ in sections \ref{section5} and \ref{section6}.

It will be shown that the variables $G(c,M)$ satisfy relations similar to
(\ref{Frecurs}), and we wish to establish for $G(c,M)$ an inequality analogous to
(\ref{Fform}), with a reversed inequality sign. We note, however, that for fixed
$c$, $B$ and $\gamma $ the inequality
\be \label{Gform1}
G(c,M)\leq Be^{M\left (1-\frac {\ln M}{c+1}
+\frac {\gamma }{c+1}\right )}
\ee
cannot hold for sufficiently large $M$, since for such $M$ the right-hand side of
(\ref{Gform1}) becomes smaller than 1, while the left-hand side of (\ref{Gform1}) is
clearly larger than 1.

We henceforth focus on the function $G(c,M)$ defined in (\ref{GcMsum}).
Our goal is to estimate the value of $G(c,M)$ for pairs $(c,M)$ which belong to
the domain
\be \label{Drangbet}
{\cal D}={\cal D}_{\beta }
\ee
for some $0<\beta <1/2$ (recall (\ref{Dbetadfn})).
We denote
\be \label{D1defntn}
{\cal D}_{+}=\{(c,M):e^{\beta c}<M<e^{\beta (c+1)}\}
\ee
and
\be \label{epsval}
\epsilon =1-2\beta .
\ee
Analogous to
(\ref{Frecurs}), for points $(c,M)\in {\cal D}$ we have the following relation
\be \label{Grecurs}
G(c,M)=\sum _{z=0}^{[M/c]}G(c-1,M-cz)\frac {m_0^z}{z!}e^{z^2/m_0}.
\ee
(Of course, even though $(c,M)\in {\cal D}$, some points $(c-1,M-cz)$ in
(\ref{Grecurs}) may fail to belong to ${\cal D}$.)

To obtain an upper bound of the type (\ref{Gform1}) on ${\cal D}$  we will employ the
iterative method described in sections \ref{section5} and \ref{section6}. To use this
approach in the present situation we have to guarantee in advance that (\ref{Gform1}) holds
for points in ${\cal D}_{+}$. This property will be a consequence of the following results.

\begin{Prop}\label{nurecrse}
Let $\{p_k\}_{k=1}^{\infty }$ denote the sequence of primes. Then
\be \label{nupkrecurs}
\psi (x,p_{k+1})=\sum _{j=0}^{N_{k+1}}\psi
\left (x/p_{k+1}^j,p_k\right )
\ee
holds for every $x>2$ and $k\geq 1$, where we denote
${\displaystyle N_k=\left [\frac {\ln x}{\ln p_{k+1}}\right ]}$.
\end{Prop}
{\em Proof}: Let ${\cal F}_k(x)$ denote the set of integers $z\leq x$ whose largest
prime divisor does not exceed $p_k$, so that
\be \label{nuFkNrl}
\psi (x,p_k)=\#\{{\cal F}_k(x)\}.
\ee
Denote by $A_j$ the set of integers $z\in {\cal F}_{k+1}(x)$ such that $p_{k+1}^j$
is the largest power of $p_{k+1}$ which divides $z$.
It is then easy to see that
\be \label{Ajdefinit}
A_j=p_{k+1}^j{\cal F}_k\left (\frac {x}{p_{k+1}^j}\right )
\ee
and
\be \label{Fk1union}
{\cal F}_{k+1}(x)=\bigcup _{j\geq 0}A_j,
\ee
a disjoint union. Equation (\ref{nupkrecurs}) follows
from (\ref{nuFkNrl}), (\ref{Ajdefinit}) and (\ref{Fk1union}). $\hfill \Box$

\begin{Prop}\label{upperbundln2}
Let $\alpha >1$ be fixed, and consider pairs $(x,y)$ such that
\be \label{nalpN2}
y=\alpha (\ln x)^2.
\ee
Then there exists a constant $C>1$ such that
\be \label{nlnNsqrbnd}
\frac {\ln \psi (x,y)}{\ln x}< \frac {1}{2}+\frac {C}{\ln y}
\ee
holds for every $x>1$, where $y$ is as in $(\ref{nalpN2})$.
\end{Prop}
The proof is displayed in the appendix.

\begin{Prop}\label{G2Fcm}
Let ${\cal D}$ be as in $(\ref{Drangbet})$.
Then there exist constants $K$ and $c_0$ such that
\be \label{GFnuinq}
G(c,M)<KF(c,M)
\ee
holds for every $c\geq c_0$.
\end{Prop}
{\em Proof}: We note that
$$
z\leq \frac {M}{c}\leq \frac {e^{\beta (c+1)}}{c} \mbox { and } m_0>\frac {e^c}{c},
$$
implying
\be \label{newexpbnd}
\frac {z^2}{m_0}<\frac {e}{ce^{\epsilon c}},
\ee
and it follows that
$$
e^{z^2/m_0}<\exp \{3e^{-\epsilon c}/c\}.
$$
We fix a constant $c_0$, and then (\ref{GFnuinq}) follows from (\ref{Frecurs}),
(\ref{Grecurs}) and (\ref{newexpbnd}) for $c\geq c_0$ , by employing induction on $c$.
$\hfill \Box$
\begin{Rem}\label{FreplG}
We will establish an upper bound for $F(c,M)$, and
then use $(\ref{GFnuinq})$ to estimate $G(c,M)$ from above. Thus we wish to establish
for $F$ an inequality of the form
\be \label{Fform2}
F(c,M)\leq B_1e^{M\left (1-\frac {\ln M}{c+1}
+\frac {\gamma }{c+1}\right )}
\ee
for some coefficient $B$ and a certain $\gamma $ $($which may depend on $c$ and $M)$,
and in view of $(\ref{GFnuinq})$ this will yield the estimate
\be \label{GdefnE}
G(c,M)\leq B\exp \left \{M\left (1-\frac {\ln M}{c+1}
+\frac {\gamma }{c+1}\right )\right \}.
\ee
\end{Rem}

The following result is a consequence of Proposition \ref{upperbundln2}.
\begin{Prop}\label{insiD1}
Let ${\cal D}_{+}$ be as in $(\ref{D1defntn})$, and let $C$ be as in Proposition
$\ref{upperbundln2}$. Then $(\ref{Fform2})$, with $B=1$ and $\gamma =C$,
holds on ${\cal D}_{+}$.
\end{Prop}

We consider (\ref{Grecurs}) as a difference equation in ${\cal D}$ satisfying boundary
upper bounds on ${\cal D}_{+}$ as expressed in Proposition \ref{insiD1}.
For a fixed $\kappa >1$ let
$$
E_{\kappa }={\cal D}\cap \{1\leq c\leq \kappa \}
$$
which is a bounded set, and it follows that for any fixed $\gamma $,
$F(\cdot ,\cdot)$ satisfies (\ref{Fform2}) on $E_{\kappa }$
for some $B>1$ (depending on $\gamma $).

Suppose that we have an upper bound for $F(\cdot ,\cdot )$ on $E_{\kappa }$,
and we consider in the left hand side of (\ref{Frecurs}) pairs $(c,M)$ which belong to
$E_{\kappa +1}\setminus E_{\kappa }$. We will next show that for such $(c,M)$
the right hand side of (\ref{Frecurs}) involves pairs $(c-1,M-cz)$ for which an
upper bound of the form (\ref{Fform2}) has been already established. We will then
use these bounds to estimate the right hand side of (\ref{Frecurs}), thus establishing
an upper bound for $F(c,M)$.
\begin{Prop}\label{belong}
If $(c,M)\in E_{\kappa +1}\setminus E_{\kappa }$ then
\be \label{cMzblong}
(c-1,M-cz)\in E_{\kappa }\cup {\cal D}_{+}
\ee
for every $0\leq z\leq M/c$.
\end{Prop}
{\em Proof}: If $(c,M)\in E_{\kappa +1}$ then
$\displaystyle {M\leq e^{\beta c}}$. Obviously
this can be written in the form
$$M\leq e^{\beta [(c-1)+1]},$$
implying that $(c-1,M)\in {\cal D}_{+}$ if
$\displaystyle{M>e^{\beta (c-1)}}$, and $(c-1,M)\in E_{\kappa }$
if $\displaystyle{M\leq e^{\beta (c-1)}}$. $\hfill \Box$

It follows from Proposition \ref{belong} that each summand $F(c-1,M-cz)$
in the right hand side of (\ref{Grecurs}) may be bounded by employing a bound of the form
(\ref{Fform2}) for $(c-1,M-cz)$.

In analogy with (\ref{mzeE}) we have that
\be \label{GmzeE}
\frac {m_0^z}{z!}<e^{\bar E},
\ee
where similarly  to (\ref{Elogaritm})
\be \label{GElogaritm}
\bar E=\left (z\ln m_0-z\ln z+z\right ).
\ee
(In (\ref{GElogaritm}) we ignore the term $\sqrt z$ in (\ref{Stirapprx}), since we consider
now an upper bound.) Substituting $m_0=(e-1)e^c/c$ in (\ref{GElogaritm}) we obtain
$$
\bar E= cz-z\ln z+z(1+\ln (e-1))-z\ln c.
$$
Let $A$ be as in (\ref{Alogaritm}), and analogous to (\ref{c1Mzbnd}) we assume that
$$
F(c-1,M-cz)\leq Be^A,
$$
so that
$$
F(c-1,M-cz)\frac {m_0^z}{z!}\leq Be^{A+\bar E}.
$$
It follows that an upper bound for $A+\bar E$ is given by the function
$H(z)$ in (\ref{maxexpres}), where the variable $a$ (recall (\ref{adefinitn}))
is replaced by $a^{\star }$ in (\ref{limita}).
We still denote this function by $H(z)$, and
analogous to (\ref{FcMHD}) we have the relation
\be \label{GcMHD}
F(c-1,M-cz)\frac {m_0^z}{z!}<Be^{ H(z)}.
\ee

As in section \ref{section5}, we should maximize the function $ H(z)$ over
$0\leq z\leq [M/c]$. But in the present situation, since we are concerned
with an upper bound, we may use the maximum of $H(z)$ over the real
interval $0\leq z\leq M/c$ and do not have to restrict to the integers in this interval.

Summarizing the above discussion we obtain, analogous to (\ref{Hzmaxlower}), the
following result.
\begin{Prop}\label{GcMBbond}
Assume that
\be \label{GcMBinq}
F(c,M)\leq Be^{M\left (1-\frac {\ln M}{c+1}+\frac {\gamma}{c+1}\right )}
\ee
for every $(c,M)\in E_{\kappa }$, for some $\gamma>C$ and $\kappa >1$. Then
\be \label{maxestimate}
\max \left \{F(c-1,M-cz)\frac {m_0^z}{z!}:0\leq z\leq \frac {M}{c}\right \}
\leq Be^{M\left (1-\frac {\ln M}{c}+\frac {\gamma +f(\gamma )}{c}\right )},
\ee
implying
\be \label{muineqult}
F(c,M)<Be^{M\left (1-\frac {\ln M}{c}+\frac {\gamma +f(\gamma )}{c}\right )+\ln (M/c)}
\ee
for every $(c,M)\in E_{\kappa +1}$.
\end{Prop}
\begin{Rem}
The term $\ln (M/c)$ appears in $(\ref{muineqult})$ since we should multiply the
maximum in $(\ref{maxestimate})$ by the number of terms which appear in the sum in
$(\ref{Frecurs})$. We may use $\ln (M/c)$ rather than $\ln ([M/c]+1)$ since there are
in $(\ref{Frecurs})$ several summands which are much smaller than the maximal term there.
\end{Rem}
In this section we use induction to establish an inequality of the type
(\ref{Fform2}), with $\gamma $ depending on $(c,M)$ as follows:
\be \label{gamcM7}
\gamma (c,M)=\bar a+\ln c+\ln \ln c-\ln ^{(2)}M-\ln ^{(3)}M
\ee
for a certain $\bar a>a^{\star }$.\\

We consider now the maximization
in the left hand side of (\ref{maxestimate}). Employing an induction hypothesis we obtain
bounds on the expressions $F(c-1,M-cz)$, using inequalities  of the form (\ref{GcMBinq})
for the pairs $(c-1,M')$, where $M'=M-cz$. In these bounds we denote
$\gamma =\gamma (c-1,M')$, using (\ref{gamcM7}).
Suppose that the {\em maximum over the bounds} is attained at $1<M_0\leq M$, and denote
$\gamma _0=\gamma (c-1,M_0)$, namely
\be \label{gama0defn}
\gamma _0=\bar a+\ln (c-1)+\ln \ln (c-1)-\ln ^{(2)}M_0-\ln ^{(3)}M_0.
\ee
Clearly the maximum over the bounds is not larger than the maximal value of
\be \label{gam0functn}
\exp \left \{(M-cz)\left [1-\frac {\ln (M-cz)}{c}+\frac {\gamma _0}{c}\right ]\right \}
\frac {m_0^z}{z!}
\ee
over $0\leq z\leq M/c$.

In view of (\ref{maxestimate}) and (\ref{muineqult}), and analogous to (\ref{Mgamampr1}),
we wish to establish
\be \label{gamafgamampr}
\frac {\gamma _0+f(\gamma _0)}{c}<\frac {\ln M}{c(c+1)}+\frac {\gamma '}{c+1},
\ee
where
\be \label{gamprdefn}
\gamma '=\bar a+\ln c+\ln \ln c-\ln ^{(2)}M-\ln ^{(3)}M.
\ee
We first address the term $f(\gamma _0)$ in (\ref{gamafgamampr}),
and recalling (\ref{fdefn}) and (\ref{gama0defn}) we have
\be \label{fgamaapr}
f(\gamma _0)=\ln \left (1+e^{a^{\star }-\bar a}
\frac {\ln M_0\ln ^{(2)}M_0}{(c-1)\ln (c-1)}\right ).
\ee
We assume now that
$$
(c,M)\in {\cal D}_{\beta },
$$
and denote in (\ref{gamcM7})
\be \label{abardefn}
\bar a=a^{\star}+\delta
\ee
for some $\delta >0$.
For small enough $\beta $, arguing as in Remark \ref{smallbeta} we have, analogous
to (\ref{5star})
$$
\frac {\ln M_0}{c-1}<(1+\epsilon)\beta.
$$
Clearly we have also $\frac {\ln \ln M_0}{\ln (c-1)}<1$, and thus,
if $\delta $ is sufficiently small, then
\be \label{fgamlnMest}
\frac {f(\gamma _0)}{c}<\frac {q\ln M}{c(c+1)} \mbox { if } c>c_0
\ee
for some $c_0$, where
\be \label{qualpha}
e^{-\delta }<q<1.
\ee
We note that $q$ in (\ref{qualpha}) may be arbitrarily close to $e^{-\delta }$ provided that
$\beta >0$ and $\delta >0$ are sufficiently small. Specifically we may choose the parameters
$\delta $ and $q$ in (\ref{abardefn}), (\ref{fgamlnMest}) and (\ref{qualpha}) as follows:
\be \label{deltqdfn}
\delta =\lambda \beta \mbox { and } q=1-\frac {\lambda \beta }{2}
\ee
where $\lambda >0$ may be arbitrarily small.

We next consider the terms $\gamma _0/c$ and $\gamma '/(c+1)$ in (\ref{gamafgamampr}). Let
$z_0$ be the point where the maximization over $z$ of (\ref{gam0functn}) is attained, and
let, as above, $M_0=M-cz_0$. We note that in this maximization, the value $\gamma _0$ is
the same for all the points $(c-1,M')$, $1<M'\leq M$. We have then
\be \label{M0definit}
M_0=M(1-t_0),
\ee
where by (\ref{toexpr})
$$
t_0=\frac {1}{1+e^{\gamma _0-a^{\star }}}<\frac {e^{-\delta }\ln M_0\ln ^{(2)}M_0}{c\ln c}.
$$
Thus
\be \label{lnt0estm}
\ln (1-t_0)>-\frac {q_1\ln M_0\ln ^{(2)}M_0}{c\ln c}
\ee
for some constant $e^{-\delta }<q_1<1$. It follows from (\ref{M0definit}) and
(\ref{lnt0estm}) that
$$
\left (1+\frac {q_1\ln ^{(2)}M_0}{c\ln c}\right )\ln M_0>\ln M,
$$
hence
$$
\ln M_0>\left (1-\frac {q_1\ln ^{(2)}M_0}{c\ln c}\right )\ln M,
$$
and we obtain
\be \label{lnlnMM0}
\ln \ln M_0>\ln \ln M-\frac {q_1\ln ^{(2)}M}{c\ln c}
\ee
for some constant $q_2$.

Using the expressions (\ref{gama0defn}) and  (\ref{gamprdefn}) it follows from
(\ref{lnlnMM0}) that
\newline
$\displaystyle{\frac {\gamma _0}{c}-\frac {\gamma '}{c+1}}$ is smaller than
$$
\frac {\bar a+\ln (c-1)+\ln \ln (c-1)}{c}
-\frac {\bar a+\ln c+\ln \ln c}{c+1}+\frac {q_2\ln ^{(2)}M}{c^2\ln c},
$$
implying that
\be \label{fracgam0gm2}
\frac {\gamma _0}{c}-\frac {\gamma '}{c+1}<
\frac {\bar a+\ln c+\ln \ln c}{c(c+1)}+\frac {q_2\ln ^{(2)}M}{c^2\ln c}
\ee
If $(c,M)$ is such that
$$
\ln c<(1-q)\ln M,
$$
then (\ref{gamafgamampr}) would follow from (\ref{fgamlnMest}) and (\ref{fracgam0gm2})
for large enough $c$. We thus consider pairs $(c,M)$ satisfying
\be \label{cMkappa}
M>c^{\nu }
\ee
for some constant $\nu >1$ such that
\be \label{kapaque}
(1-q)\nu \geq 1.
\ee
If we choose, as in (\ref{deltqdfn}), $q=1-\lambda \beta /2$ for some $\lambda >0$, we may take
\be \label{spcilnu}
\nu =\frac {2}{\lambda \beta }.
\ee
We have thus established the following result.
\begin{Prop}\label{discsummry}
Let $\bar a$ and $\delta >0$ be as in $(\ref{abardefn})$, let
$\gamma (c,M)$ be as in $(\ref{gamcM7})$, and consider pairs
$(c,M)\in {\cal D}$ which satisfy $(\ref{cMkappa})$ and $(\ref{kapaque})$.
Then there exist constants $B$, $c_0$ and $\delta _0$ such that
\be \label{acumineqlt}
F(c,M)<Be^{M\left (1-\frac {\ln M}{c+1}+\frac {\gamma (c,M)}{c+1}\right )+c\ln (M/c)}
\ee
holds provided that $c>c_0$ and $\delta >\delta _0$.
\end{Prop}
{\em Proof}: The inequality (\ref{acumineqlt}) follows from (\ref{muineqult}) and
(\ref{gamafgamampr}) and the preceding discussion. We note that when employing successively
the inequalities (\ref{muineqult}) and (\ref{gamafgamampr}), the various terms
$\ln (M/c)$ in (\ref{muineqult}) accumulate, yielding the term $c\ln M$ in
(\ref{acumineqlt}). $\hfill \Box$

Concerning $G(c,M)$, in view of Remark \ref{FreplG} we obtain the following result:
\begin{Thm} \label{mainuprbndrslt}
Consider pairs $(c,M)$ satisfying
\be \label{uprcondi}
c^{\nu }<M<e^{\beta c}
\ee
for some $\nu >2$ and $0<\beta <1/2$. Then there exist constants $\overline{a}>a^{\star }$
and $c_0$ such that
\be \label{nineqult3}
G(c,M)<e^{M\left (1-\frac {\ln M+\ln ^{(2)}M+\ln ^{(3)}M}{c}
+\frac {\overline{a}+\ln c+\ln ^{(2)}c}
{c}\right )}
\ee
for every $c>c_0$. Moreover, for every $\lambda >0$, which may be arbitrarily small, we may take
$$
\overline{a}=a^{\star }+\lambda \beta
$$
provided that $\beta >0$ is sufficiently small and $\nu \geq 2/\lambda \beta $.
\end{Thm}
The last assertion of the theorem follows from (\ref{spcilnu}).
\mysection{The main results}\label{section8}
In this section we will establish our main results concerning lower and upper bounds
for $\psi (x,y)$. They consist of rephrasing the results in sections \ref{section6}
and \ref{section7} in terms of $x$ and $y$ instead of $c$ and $M$. We obtain from Theorem
\ref{mainlwrbndrslt} our first main result:
\begin{Thm}\label{noaymptot8}
Consider $(x,y)$ such that
\be \label{lowrange}
\exp \{(\ln y)^{1-\theta }\}<\ln x<\sqrt {y}
\ee
for some $\theta >0$. Then there exists an $\underline{a}$ and a $y_0$ such that
\be \label{lwrboundassy8}
\frac {\ln \psi (x,y)}{\ln x}>1-\frac {\ln ^{(2)}x+\ln ^{(3)}x+\ln ^{(4)}x}{\ln y}
+\frac {\underline{a}+\ln ^{(2)}y+\ln ^{(3)}y}{\ln y}
\ee
for every $(x,y)$ satisfying $(\ref{lowrange})$ and $y>y_0$.
\end{Thm}

Concerning an upper bound for $\psi (x,y)$, Theorem \ref{mainuprbndrslt} yields our second
main result:

\begin{Thm}\label{discsummry}
For some constants $0<\beta <1/2$ and $\nu >0$ consider pairs $(x,y)$ which satisfy
\be \label{upperange}
(\ln y)^{\nu }<\ln x<y^{\beta },
\ee
and let $u$ be as in $(\ref{unotatn})$. Then there exist constants $y_0$ and
$\overline{a}>a^{\star }$ such that
\be \label{acumineqlt8}
\frac {\ln \psi (x,y)}{\ln x}<1-\frac {\ln ^{(2)}x+\ln ^{(3)}x+\ln ^{(4)}x}{\ln y}
+\frac {\overline{a}+\ln ^{(2)}y+\ln ^{(3)}y}{\ln y}+\frac {\ln u}{u}
\ee
holds provided that $y>y_0$. If in $(\ref{upperange})$ we have $\nu >2 $ then
\be \label{acumineqlt18}
\frac {\ln \psi (x,y)}{\ln x}<1-\frac {\ln ^{(2)}x+\ln ^{(3)}x+\ln ^{(4)}x}{\ln y}
+\frac {\overline{a}+\ln ^{(2)}y+\ln ^{(3)}y}{\ln y}
\ee
holds for every $y>y_1$, for some $y_1$. Moreover, we may take $\overline{a}>a^{ \star }$
to be arbitrarily close to $a^{ \star }$ provided that $\beta $ is small enough and
$\nu $ is large enough.
\end{Thm}

\begin{Rem}\label{conject}
The bounds $(\ref{lwrboundassy8})$ and $(\ref{acumineqlt18})$ raise the conjecture that
for each $k\geq 2$, in a certain range of the variables $x$ and $y$ the following bounds
\be \label{acumineqlt81}
\frac {\ln \psi (x,y)}{\ln x}>1-\frac {1}{\ln y}\left [
\sum _{j=2}^{k+1}\ln ^{(j)}x-\underline{a}-\sum _{j=2}^k\ln ^{(j)}y \right ]
\ee
and
\be \label{acumineqlt82}
\frac {\ln \psi (x,y)}{\ln x}<1-\frac {1}{\ln y}\left [
\sum _{j=2}^{k+1}\ln ^{(j)}x-\overline{a}-\sum _{j=2}^k\ln ^{(j)}y \right ]
\ee
are valid for certain constants $\underline{a}$ and $\overline{a}$.
\end{Rem}
\begin{Rem}\label{newform}
The inequalities $(\ref{lwrboundassy8})$ and
$(\ref{acumineqlt18})$ may be written in the form
\be \label{newlwrboundassy8}
\ln \left (\frac {\psi (x,y)}{x}\right )>-u[\ln u+\ln ^{(3)}x
-\ln ^{(3)}y+\ln ^{(4)}x-\underline{a}]
\ee
and
\be \label{newacumineqlt18}
\ln \left (\frac {\psi (x,y)}{x}\right )<-u[\ln u+\ln ^{(3)}x
-\ln ^{(3)}y+\ln ^{(4)}x-\overline{a}]
\ee
respectively, where we used
\be \label{ln2xln2ylnu}
\ln ^{(2)}x-\ln ^{(2)}y=\ln u.
\ee
\end{Rem}

We will next estimate the value of iterated logarithms $\ln ^{(k)}x$ and $\ln ^{(k)}y$
for pairs $(x,y)$ which satisfy
\be \label{kapbetrang}
\exp (\ln y)^{\nu }<\ln x<y^{\beta }
\ee
for some $0<\nu <1$.
To do this we will use the iterated logarithms $\ln ^{(k)}u$.
\begin{Prop}\label{iterlogest}
Let $(x,y)$ be such that $(\ref{kapbetrang})$ holds, and let $u$ be as in $(\ref{unotatn})$.
Then
\bigskip

\noindent
$(i)$ For every $k\geq 3$
\be \label{logxiter}
\ln ^{(k)}x=\ln ^{(k-1)}u+o(1).
\ee
\bigskip

\noindent
$(ii)$ For every $k\geq 4$
\be \label{logyiter}
\ln ^{(k)}y=\ln ^{(k)}u+o(1),
\ee
and
\be \label{logy3iter}
\ln ^{(3)}u+o(1)<\ln ^{(3)}y<\ln ^{(3)}u+\ln (1/\nu )+o(1).
\ee
\end{Prop}
{\em Proof}: It follows from the left inequality in (\ref{kapbetrang}) that
\be \label{lny2bnd}
\ln ^{(2)}y<\frac {1}{\nu }\ln ^{(3)}x.
\ee
We conclude from (\ref{ln2xln2ylnu}) and (\ref{lny2bnd}) that
\be \label{lnx2lnu}
\ln u<\ln ^{(2)}x<\ln u+\frac {1}{\nu }\ln ^{(3)}x.
\ee
Since
$$\frac {\ln ^{(3)}x}{\ln ^{(2)}x}=o(1) \mbox { as } x\to \infty,$$
it follows from (\ref{lnx2lnu}) that
$$
\ln ^{(2)}x=(\ln u)(1+o(1)),
$$
which establishes (\ref{logxiter}) for every $k\geq 3$.

Concerning $\ln ^{(3)}y$ we obtain from the right inequality in (\ref{kapbetrang}) that
\be \label{betlnylnu}
\beta \ln y>\ln u+\ln ^{(2)}y.
\ee
Since
$$
\frac {\ln ^{(2)}y}{\ln y}=o(1) \mbox { as } y\to \infty,
$$
we conclude from (\ref{betlnylnu}) that
$$
\beta (1+o(1))\ln y>\ln u,
$$
and  therefore
$$
\ln ^{(2)}y>\ln ^{(2)}u-\ln \beta +o(1),
$$
implying
\be \label{frstpart}
\ln ^{(3)}y>\ln ^{(3)}u+o(1).
\ee
Moreover, it follows from from the left inequality in (\ref{kapbetrang}) that
\be \label{uln2ylnk}
\ln u+\ln ^{(2)}y>(\ln y)^{\nu }.
\ee
Since
$$
\frac {\ln ^{(2)}y}{(\ln y)^{\nu }}=o(1)
$$
we obtain
$$
(\ln y)^{\nu }<(1+o(1))\ln u.
$$
This implies
$$
\ln ^{(3)}y<\ln ^{(3)}u-\ln \nu +o(1),
$$
which together with (\ref{frstpart}) establishes (\ref{logy3iter}). The relations
(\ref{logyiter}) for $k\geq 4$ follow from (\ref{logy3iter}). The proof of the
proposition is complete. $\hfill \Box$

The estimates in Proposition \ref{iterlogest} yield the approximation of
$\ln \rho (u)$ presented in Corollary \ref{lnroest}.

We conclude this section by employing Theorems \ref{noaymptot8} and \ref{discsummry} to
establish a result concerning Bertrand's Conjecture. As is well known, Bertrand's conjecture
was that for every integer $y$ there exists a prime $p$ satisfying $y\leq p\leq 2y$.
\begin{Cor}\label{Bertran}
Let $\gamma >3/2$ be fixed. Then there exists a $y_0$ such that for every integer $y>y_0$
there exists a prime $p$ satisfying
\be \label{betraninq}
y<p< \gamma y.
\ee
\end{Cor}
{\em Proof}: By Remark
\ref{alpbetrel} and Theorem \ref{mainuprbndrslt} we may assume that
$$|
a^{\star}-\underline{a}|<\frac {\beta }{2} \mbox { and }
|a^{\star}-\overline{a}|<\left (\gamma -\frac {3}{2}\right )\beta
$$
provided that $\beta $ is small enough,
and that $\nu $ and $\ln x/\ln y$ are large enough. We thus assume that
the latter parameters were chosen such that
\be \label{adiffrnc}
|\overline{a}-\underline{a}|<(\gamma -1)\beta,
\ee
and such that there exists a $y_0$ for which
\be \label{cordomain}
\exp \{(\ln y)^{\nu }\}<\ln x<y^{\beta }
\ee
holds for $y=y_0$ some $x_0$. Then (\ref{cordomain}) holds for every $y>y_0$, and we
may assume that $y_0$ and $x_0$ were chosen such that $\ln x_0/\ln y_0$ is sufficiently
large, as required. For $y_1>y_0$ denote $y_2=\gamma y_1$, and let $x$ be such that
both $(x,y_1)$ and $(x,y_2)$ satisfy (\ref{cordomain}).

We write the inequalities (\ref{lwrboundassy11n}) and (\ref{nineqult31n}) in the form
\be \label{lnpsiform}
\ln \psi (x,y)=u[\ln y+\ln ^{(2)}y+\ln ^{(3)}y-\ln ^{(2)}x-\ln ^{(3)}x-\ln ^{(4)}x+a]
\ee
where $a$ satisfies $\underline{a}<a<\overline{a}$, and employ (\ref{lnpsiform}) to estimate
$\psi (x,y_2)-\psi (x,y_1)$. For a fixed value of $a$ we denote by $\psi _a(x,y)$ the
expression for $\psi (x,y)$ in (\ref{lnpsiform}). To estimate $\psi _a(x,y_2)-\psi _a(x,y_1)$
we consider  the partial derivative $(\ln \psi _a(x,y))_y$, which is equal to
\be \label{derivnpsi}
\ln x\frac {\partial }{\partial y}\left (\frac {1}{\ln y}
[\ln y+\ln ^{(2)}y+\ln ^{(3)}y-\ln ^{(2)}x-\ln ^{(3)}x-\ln ^{(4)}x+a]\right ).
\ee
It is easy to see that the expression (\ref{derivnpsi}) is larger than
\be \label{lnpsidrvtv}
\frac {\ln x\ln u}{y_1(\ln y_1)^2} \mbox { for every }y_1\leq y\leq y_2.
\ee
The fact that $(\ln \psi _a(x,y))_y$ is larger than the expression in  (\ref{lnpsidrvtv}) 
implies that
\be \label{psiaestim}
\psi _a(x,y_2)>\psi _a(x,y_1)\exp\left \{\frac {(\gamma -1)u\ln u}{\ln y_1}\right \},
\ee
where we used $y_2-y_1=(\gamma -1)y_1$

Returning to (\ref{lnpsiform}) let $a_1$ and $a_2$ correspond to $y_1$ and $y_2$ in
this formula, so that by (\ref{adiffrnc})
$$
|a_2-a_1|<
|\overline{a}-\underline{a}|<(\gamma -1)\beta,
$$
and we write
\be \label{a1a2diff}
|a_2-a_1|=\sigma \beta ,\;\sigma <\gamma -1.
\ee
It follows from (\ref{lnpsiform}), (\ref{psiaestim}) and (\ref{a1a2diff}) that
\be \label{psiaestim}
\psi (x,y_2)>\psi (x,y_1)\exp\left \{\frac {(\gamma -1)u\ln u}{\ln y_1}
-\sigma \beta u
\right \}.
\ee
By (\ref{cordomain}) the pair $(x,y_1)$ satisfies
$$
\frac {\ln \ln x}{\ln y_1}<\beta,
$$
and moreover, taking $x$ sufficiently large we can have $\ln ^{(2)}x/\ln y_1$ be arbitrarily
close to $\beta $. In this case we also have
\be \label{lnulny1}
\left |1-\frac {1}{\beta }\frac {\ln u}{\ln y_1}\right | \mbox { is arbitrarily small}
\ee
provided that $x$ is sufficiently large.
Writing the exponent in the right hand side of (\ref{psiaestim})
in the form
\be \label{exponform}
(\gamma -1)\beta u\left [\frac {1}{\beta }\frac {\ln u}{\ln y_1}
-\frac {\sigma }{\gamma -1}\right ]
\ee
yields, in view of (\ref{lnulny1}) and $\sigma <\gamma -1$, that
$$
\psi (x,y_2)>2\psi (x,y_1),
$$
if $x$ is large enough, from which we conclude that
\be \label{psi12dffrnc}
\psi (x,y_2)-\psi (x,y_1)>2.
\ee
But clearly (\ref{psi12dffrnc} implies that there exists a prime $p$ satisfying
$y_1<p<y_2$. This establishes (\ref{betraninq}), and completes the proof of the corollary.
$\hfill \Box$.

\mysection{Appendix}\label{section9}
{\em Proof of Theorem $\ref{niceresult}$}: Let $F=[1,x]\setminus E$ be the complement of $E$
in $[1,x]$. For a prime
$\sqrt x\leq p\leq x$ we denote by $F_p$ the set of integers in $F$ which are divisible
by $p$. Then $F_{p_1}\cap F_{p_2}=\emptyset $ if $p_1\not =p_2$,
$$
\#(F_p)=\left [\frac {N}{p}\right ]
$$
and it follows that
\be \label{Npsum}
\#(F)=\sum _{\sqrt x\leq p\leq x}[x/p]<x\sum _{p\geq \sqrt x}^x\frac {1}{p},
\ee
where the sum is over the primes in the indicated interval.
To estimate the sum in the right hand side of (\ref{Npsum}) we consider, more generally,
sums of the form
\be \label{anbsum}
S_{a,b}=\sum _{a\leq p\leq b}\frac {1}{p}.
\ee
By the Prime Numbers Theorem the distribution function of the number of primes in
the real line is, for large
enough $t$, $\Phi (t)=t/\ln t$.
Using this in the summation in (\ref{anbsum}) implies that for sufficiently large $a$
we have
$$
S_{a,b}\approx \int _a^b\frac {d\Phi (t)}{t}= \int _a^b\frac {\Phi(t)dt}{t^2}
+\left.\frac {\Phi (t)}{t}\right |_a^b,
$$
and substituting $\Phi (t)=t/\ln t$ we conclude that
\be \label{Sabsum}
S_{a,b}\approx \int _a^b\frac {dt}{t\ln t}
+\left.\frac {1}{\ln t}\right |_a^b <\ln \ln b-\ln \ln a.
\ee
For $a=\sqrt x$ and $b=x$ the right hand side of (\ref{Sabsum}) is equal to $\ln 2$,
and using this in (\ref{Npsum}) yields that for sufficiently large $x$ we have
$$\#(F)<x\ln 2,$$ implying $$\#(E)>x\ln (e/2).$$ This establishes (\ref{FsqrtN}) and
concludes the proof. $\hfill \Box$.

{\em Proof of Proposition} \ref{nlnNsqrbnd}:
It follows from $\psi (x,2)\leq \ln x/\ln 2$ that
$$
\psi (x,2)\leq \frac {\sqrt x}{\ln 2},
$$
since $\ln x<\sqrt x$ for every $x\geq 1$.
It is easy to see that
\be \label{nuprod}
\psi (x,p_k)\leq \frac {\sqrt x}{(\ln 2)(1-1/\sqrt {p_2})\cdots (1-1/\sqrt {p_k})},
\ee
for every $k\geq 2$. Relation (\ref{nuprod}) can be established by
employing a simple induction argument, using (\ref{nupkrecurs}).

To estimate from above the right hand side of (\ref{nuprod}), we have to estimate
from below the product
\be \label{prodestim}
\prod _{j=1}^k\left (1-\frac {1}{\sqrt {p_j}}\right ),
\ee
and for this we estimate from above the sum
\be \label{sumestim}
\sum _{j=1}^k\frac {1}{\sqrt {p_j}}.
\ee
To this end we use the distribution function
$$
\Phi (t)=\frac {t}{\ln t}
$$
of the primes in the real line, and we have to estimate
$$
\int _3^{p_k }\frac {d\Phi (t)}{\sqrt t}.
$$
This leads to
\be \label{CpklnlnN}
\int _3^{p_k }\frac {dt}{\sqrt t\ln t}=\int _{\sqrt 3}^{\sqrt {p_k}}\frac {ds}{2\ln s}<
\frac {C\sqrt {p_k}}{\ln p_k}
\ee
for some constant $C>0$, and we obtain
\be \label{nupkNsq}
\psi (x,p_k)\leq x^{1/2}\exp\left\{C\sqrt {p_k}/\ln p_k\right \}.
\ee
For a prescribed $y=\alpha (\ln x)^2$
we let $p_k$ be the smallest prime $p$ which satisfies $p\geq y$. Employing
(\ref{nupkNsq}) for this $p_k$ yields
the assertion of the proposition.
$\hfill \Box$


\begin{thebibliography}{9999}
\bibitem{Bret} de la Bret\`{e}che R. and G. Tenenbaum (2002). Local distribution
of the $k$th divisor of an integer, {\em Proc. London Math. Soc.}, {\bf 2}, 289-323.


\bibitem{Bruij} de Bruijn N. G. (1951). On the number of positive integers $\leq x$ and free
of prime factors $>y$, {\em Nederl. Akad. Wetensch. Proc. Ser. A}, {\bf 54}, 50-60.

\bibitem{Dick} Dickman K. (1930). On the frequency of numbers containing prime factors
of a certain relative magnitude, {\em Ark. Mat. Astr. Fys.}, {\bf 22}, 1-14.

\bibitem{Erd1} Erd\H{o}s P. (1935). On the normal number of prime factors of $p-1$ and some other
related problems cocerning Euler's $\Phi $-function, {\em Quart. J. Math. (Oxford)},
{\bf 6}, 205-213.

\bibitem{Erd2} Erd\H{o}s P. (1952). On the greatest prime factor of $\prod f(k)$,
{\em J. London Math. Soc.}, {\bf 27}, 379-384.

\bibitem{Erd3} Erd\H{o}s P. (1955). On consecutive integers, {\em Nieuw Arch. Wisk.} (3),
{\bf 3}, 124-128.

\bibitem{ErdS} Erd\H{o}s P. and A. Schinzel (1990). On the greatest prime factor of
$\prod f(k)$, {\em Acta Arith.}, {\bf 55}, 191-200.

\bibitem{Fouv} Fouvry E. and G. Tenenbaum (1996). Statstical distribution of integers
without prime factors in arithmetic progressions, {\em Proc. London Math. Soc.}
{\bf 3}, 481-514.

\bibitem{Frid1} Friedlander J. B. (1973). Integers without large prime factors, {\em Nederl.
Wetensch. Proc. Ser. A}, {\bf 76}, 443-451.

\bibitem{Frid2} Friedlander J. B. (1976). Integers free from large and small primes,
{\em Proc. London Math. Soc.} (3), {\bf 33}, 565-576.

\bibitem{Frid3} Friedlander J. B. (1981). Integers without large prime factors II,
{\em Acta Arith.}, {\bf 39}, 53-57.

\bibitem{Frid4} Friedlander J. B. (1984). Integers without large prime factors III,
{\em Arch. Math.}, {\bf 43}, 32-36.

\bibitem{Gran1} Granville A. (1989). On positive integers $\leq x$ with prime factors
$\leq t\log x$, in {\em Number Theory and Applications} (R. A. Mollin, ed.), Kluver,
403-422.

\bibitem{Gran2} Granville A. (1991). On pairs of coprime integers with no large prime
factors, {\em Expos. Math.}, {\bf 9}, 335-350.

\bibitem{Gran3} Granville A. (1993). Integers without large prime factors in arithmetic
progressions I, {\em Acta Math.}, {\bf 170}, 255-273.

\bibitem{Gran4} Granville A. (1993). Integers without large prime factors in arithmetic
progressions II, {\em Philos. Trans. Roy. Soc. London Ser. A}, {\bf 1676}, 349-362.

\bibitem{Hazle} Hazlewood D. G. (1973). On integers all of whose prime factors are small,
{\em Bull. London Math. Soc.} {\bf 5}, 159-163.

\bibitem{Hild1} Hildebrand A. (1984). Integers free of large prime factors and the Riemann
Hypothesis, {\em Mathematika}, {\bf 31}, 258-271.

\bibitem{Hild2} Hildebrand A. (1985). Integers free of large prime factors in short
intervals, {\em Quart. J. Math. (Oxford)} (2), {\bf 36}, 57-69.

\bibitem{Hild3} Hildebrand A. (1986a). On the number of positive integers $\leq x$
and free of prime factors $>y$, {\em J. Numer Theory}, {\bf 22}, 289-307.

\bibitem{Hild4} Hildebrand A. (1986b). On the local behavior of $\Psi (x,y)$, {\em Trans.
Amer. Math. Soc.}, {\bf 297}, 729-751.

\bibitem{Hild5} Hildebrand A. (1987). On the number of prime factors of integers without
large prime divisors, {\em J. Number Theory}, {\bf 25}, 81-106.

\bibitem{HildTen1} Hildebrand A. and G. Tenenbaum (1986). On integers free of large prime
factors, {\em Trans. Amer. Math. Soc.}, {\bf 296}, 265-290.

\bibitem{HildTen2} Hildebrand A. and G. Tenenbaum (1993). On a class of difference
differential equations arising in nuber theory, {\em J. d'Analyse Math.}, {\bf 61},
145-179.

\bibitem{HildTen3} Hildebrand A. and G. Tenenbaum (1993). Integers without large
prime factors, {\em J. de Th\'{e}orie des Nombres de Bordeaux}, {\bf 2}, 411-484.

\bibitem{Hunt} Hunter S. and J. Sorenson (1997). Approximating the number of integers
free of large prime factors.
{\em Math. Comp.}, {\bf 220}, AMS, 1729-1741.

\bibitem{Huxl} Huxley M. (1973). The difference between consecutive primes.
{\em Proceeding of the Symposium in Pure athematics}, {\bf 24}, AMS, 141-146.

\bibitem{Pomer1} Pomerance C. (1980). Popular values of Euler's function,
{\em Matematika}, {\bf 27}, 84-89.

\bibitem{Pomer2} Pomerance C. (1987). Fast, rigorous factorization and discrete logrithm
algorithms, in: {\em Discrete Algorithms and Complexing (Kyoto, 1986)}, Academic Press,
Boston.

\bibitem{Ramac1} Ramachandra K. (1969). A note on numbers with a large prime factors,
{\em J. London Math. Soc.} (2), {\bf 1}, 303-306.

\bibitem{Ramac2} Ramachandra K. (1970). A note on numbers with a large prime factors II,
{\em J. Indian Math. Soc. (N. S)}, {\bf 34}, 39-48.

\bibitem{Ramac3} Ramachandra K. (1970). A note on numbers with a large prime factors III,
{\em Acta Arith.}, {\bf 19}, 49-62.

\bibitem{Rank} Rankin R. A. (1938). The difference between consecutive prime numers,
{\em J. London Math. Soc.}, {\bf 13}, 242-247.

\bibitem{Song} Song J. M. (2002). Sums of nonnegative multiplicative functions over
integers without large prime factors II. {\em Acta. Arith.}, {\bf 2}, 105-129.

\bibitem{Suzuk1} Suzuki K. (2004). An estimate for the number of integers without large
factors, {\em Math. Comp.}, {\bf 246}, 1013-1022.

\bibitem{Suzuk2} Suzuki K. (2006). Approximating the number of integers without large
prime numbers, {\em Math. Comp.}, {\bf 254}, 1015-1024.

\bibitem{Scour} Scourfiled E. J. (2006). On ideals free of large prime factors,
{\em J. de Th\'{e}orie des Nombres de Bordeaux}, {\bf 3}, 733-772.

\bibitem{Tenen1} Tenenbaum G. (1985). Sur les entries sans grand facteur premier, in:
S\'{e}minair de Th\'{e}orie des Nombres, Bordeaux 1984-85, Univ. Bordeaux I, Talence.

\bibitem{Tenen2} Tenenbaum G. (1990). Sur un probl\`{e}me d'Erd\H{o}s et. Alladi, in:
S\'{e}minair de Th\'{e}orie des Nombres (C. Goldstein ed.), Paris 1988-89,
Birkhauser, {\em Progress in Math.}, {\bf 91}, 221-239.

\bibitem{Tenen3} Tenenbaum G. (2000). A rate estimate in Billingsley's theorem for
the size distribution of large prime factors, {\em Q. J. Math.}, {\bf 3}, 385-403.

\bibitem{Versh} Vershik A. M. (1987). The asymptotic disribution of factorization of
natural numbers into prime divisors, {\em Soviet math. Dokl.}, {\bf 34}, 57-61.

\bibitem{Xu1} Xuan T. Z. (1991). The average of the divisor function over integers without
large prime factors, {\em Chinese Ann. of Math. Ser. A}, {\bf 12}, 28-33.

\bibitem{Xu2} Xuan T. Z. (1993). On the asymptotic behavior of the Dickman-de Bruijn
function, {\em Math. Ann.}, {\bf 297}, 519-533.












\end{thebibliography}
\end{document}